\documentclass[10pt, oneside]{ucthesis}
\usepackage{amsfonts,amssymb,amsmath,amsthm}
\usepackage{cite}
\usepackage{dptsymb}
\usepackage[dvips, final]{graphics}
\usepackage{isolatin1}
\def\draft{n}
\newcommand{\RR}{\mathbb R}
\newcommand{\CC}{\mathbb C}
\newcommand{\ZZ}{\mathbb Z}

\newcommand{\comma}{\mathbin ,}
\newcommand{\sltwo}{{{\mathfrak sl}_2}}
\newcommand{\frakg}{{\mathfrak g}}

\renewcommand{\qed}{~\hfill$\square$}

\newcommand{\isom}{\simeq}
\newcommand{\nisom}{\not\simeq}
\newcommand{\superset}{\supset}
\newcommand{\bigcircle}{\bigcirc}
\newcommand{\contains}{\ni}
\newcommand{\tensor}{\otimes}

\newcommand{\lbracket}{[}
\newcommand{\rbracket}{]}

\newcommand{\ad}{\operatorname{ad}}

\newcommand{\End}{\operatorname{End}}

\newcommand{\Sym}{\operatorname{Sym}}
\newcommand{\tr}{\operatorname{tr}}
\renewcommand{\det}{\operatorname{det}\nolimits}
\DeclareMathOperator{\Hom}{Hom}

\newcommand\strutn[2]{{\,^{#1}\!\!\frown^{#2}}}
\newcommand\strutu[2]{{\,_{#1}\!\!\smile_{#2}}}
\newcommand\strutv[2]{{\mathop\mid^{#1}_{#2}}}
\newcommand{\Strutn}{{\!\frown}}

\newcommand{\connect}{\mathbin \#}

\newcommand{\cA}{{\mathcal A}}
\newcommand{\cB}{{\mathcal B}}
\newcommand{\cK}{{\mathcal K}}

\newcommand{\udot}{\mathbin{\mathaccent\cdot\cup}}
\newcommand{\directedcircle}{{\circlearrowright}}

\theoremstyle{plain}
\newtheorem{theorem}{Theorem}
\newtheorem{proposition}{Proposition}[chapter]
\newtheorem{lemma}[proposition]{Lemma}
\newtheorem{corollary}[proposition]{Corollary}

\theoremstyle{definition}
\newtheorem{definition}[proposition]{Definition}
\newtheorem{relation}{Relation}

\theoremstyle{remark}

\newtheorem{exercise}[proposition]{Exercise}
\newtheorem{hint}[proposition]{Hint}
\newtheorem{remark}[proposition]{Remark}

\DeclareMathOperator{\Diag}{Diag}
\DeclareMathOperator{\Vac}{(Vacuum)}
\DeclareMathOperator{\Act}{Act}

\DeclareMathOperator{\DCable}{DCable}
\DeclareMathOperator{\CCable}{CCable}

\newcommand{\Zero}{0}
\newcommand{\One}{1}
\newcommand{\Two}{2}
\newcommand{\hopf}[2]{{\,_{#1}\!\Hopf_{#2}}}
\newcommand{\Hopf}{\HopfLink}
\newcommand{\openhopf}[2]{\OpenHopf_{\!#1}^{#2}}

\newcommand{\bc}{{\text{bc}}}

\newcommand{\bigp}[1]{\left(#1\right)}

\def\lbl#1{\label{#1}}

\def\printname#1{
	\if\draft y
		\smash{\makebox[0pt]{\hspace{-0.5in}
			\raisebox{8pt}{\tt\tiny #1}}}
	\fi
}

\newcommand{\eepic}[2]{}
\newcommand{\silenteepic}[2]{}

\newread\testin
\def\maybeinput#1{
\openin\testin=#1
\ifeof\testin\typeout{Warning: input #1 not found}\else\input#1\fi
\closein\testin
}

\newwrite\transout
\openout\transout=figs.list
\newcommand{\figscaled}[2]{%
\write\transout{-m #2 draws/#1.fig}%
\maybeinput{draws/#1.tex}%
}
\newcommand{\fig}[2]{\figscaled{#1}{#2}}
\newcommand{\figcent}[2]{\mathcenter{\figscaled{#1}{#2}}}

\usepackage{calc}
\def\mathcenter#1{
	\raisebox{.5ex+\depth-.5\totalheight}{\hbox{#1}}
}

\newcommand{\mlarge}[1]{{\text{\Large $#1$}}}

\newcounter{temp}

\usepackage{subfigure}

\begin{document}
\title{Wheeling: A Diagrammatic Analogue of the Duflo Isomorphism}
\author{Dylan P. Thurston}
\degreeyear{2000}
\degreesemester{Spring}
\degree{Doctorate of Philosophy}
\prevdegrees{A.B., Harvard University}
\chair{Vaughan Jones}
\othermembers{Robion C Kirby \\
Greg Kuperberg \\
Bruno Zumino}
\numberofmembers{4}
\field{Mathematics}
\campus{Berkeley}
\date{April 19, 2000}
\maketitle

\ssp
\renewcommand{\dsp}{}

\abstract
\dsp
We construct and prove a diagrammatic version of the Duflo isomorphism
between the invariant subalgebra of the symmetric algebra of a Lie
algebra and the center of the universal enveloping algebra.  This
version implies the original for metrized Lie algebras (Lie algebras
with an invariant non-degenerate bilinear form).  As an application of
this isomorphism, we will compute the Kontsevich integral of the
unknot and the Hopf link to all orders.

At the core of the proof, we use an elementary property of the Hopf
link which can be summarized by the equation ``$1+1=2$'' in abacus
arithmetic: doubling one component of the Hopf link is equivalent to
taking the connected sum of two Hopf links.  This property of the Hopf
link turns out, when suitably interpreted, to be exactly the property
required for the Duflo map to be multiplicative.

To compute the Kontsevich integral of the unknot, we use a property of 
the unknot that can be summarized by ``$n\cdot 0 = 0$'': the $n$-fold
connected cabling of the unknot is again an unknot.
\abstractsignature
\endabstract

\begin{frontmatter}
\tableofcontents

\begin{acknowledgements} The proof of the Wheeling theorem (the Duflo
isomorphism) was joint work with Dror Bar-Natan.  The proof of the
Wheels theorem (the Kontsevich integral of the unknot) was joint work
with Thang T.~Q.~Le.  I would like to thank these two collaborators,
with whom it was a real pleasure to work.  Many thanks also to Michel
Duflo, Stavros Garoufalidis, Greg
Kuperberg, Lev Rozansky, Michèle
Vergne, and Pierre Vogel for many
useful discussions, and especially to my advisor, Vaughan Jones.

This work was supported in part by an NSF Graduate Research Fellowship, a
Sloan Foundation Dissertation-Year Fellowship, and the Swiss National
Science Foundation.
\end{acknowledgements}
\end{frontmatter}

\dsp 
\chapter{Introduction}

\section{Elementary knot theory}\lbl{sec:ElemKnotTheory}
We begin by recalling two facts from elementary knot theory.  These
simple statements have deep consequences for Lie algebras and
Vassiliev invariants.  The two facts can be summarized by the catch
phrases ``$1+1=2$'' and ``$n\cdot 0 = 0$.''

\begin{itemize}
\item ``$1+1=2$.'' This refers to a fact in ``abacus
arithmetic.''  On an abacus, the number $1$ is naturally represented
by a single bead on a wire, as in Figure~\ref{fig:One}, which we think
of as a tangle.  The fact that $1+1 = 2$ then becomes the equality of
the two tangles in Figure~\ref{fig:OneOneTwo}.  On the left side of the
figure, ``$1+1$'', the two beads are well separated, as for connect
sum of links or multiplication of tangles; on the right side, ``$2$'',
we instead start with a single bead and double it, so the two beads
are very close together.

In other terms, the connected sum of two Hopf links is the same as
doubling one component of a single Hopf link, as in
Figure~\ref{fig:HopfSum}.
\begin{figure}[hbtp]
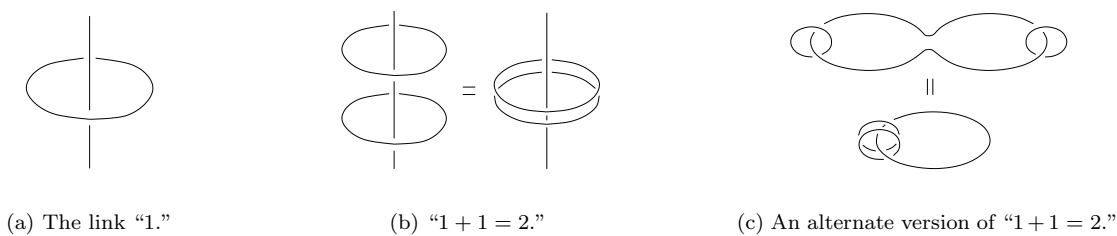
%
\begin{center}%
  \subfigure[The link ``$\One$.'']%
   {\makebox[.22\linewidth]{\figscaled{One}{0.4}\lbl{fig:One}}}%
\hfil\subfigure[``$\One + \One = \Two$.'']%
   {\makebox[.35\linewidth]{\figscaled{OneOneTwo}{0.33}\lbl{fig:OneOneTwo}}}%
\hfil\subfigure[An alternate version of ``$\One + \One = \Two$.'']%
   {\makebox[.35\linewidth]{\figscaled{HopfSum}{0.36}\lbl{fig:HopfSum}}}%
\end{center}%
\caption{Elementary knot theory, part 1}\lbl{fig:Abacus}%
\end{figure}

\item ``$n\cdot 0 = 0.$''  In the spirit of abacus
arithmetic, $0$ is represented as just a single vertical strand.  We 
prefer to close it off, yielding the knot in Figure~\ref{fig:Zero}.
The knot $n\cdot 0$ is then this knot repeated $n$ times, as in
Figure~\ref{fig:nZero}.%
The two knots are clearly the same, up to framing.

\begin{figure}[hbtp]
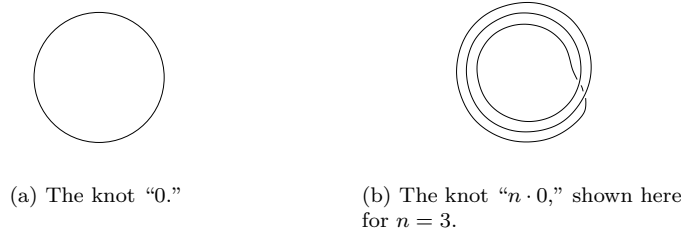
%
\begin{center}%
  \subfigure[The knot ``$\Zero$.'']%
   {\makebox[.3\linewidth]{\figscaled{Zero}{0.25}\lbl{fig:Zero}}}%
\qquad\subfigure[The knot ``$n\cdot\Zero$,'' shown here for $n=3$.]%
   {\makebox[.3\linewidth]{\figscaled{nZero}{0.25}\lbl{fig:nZero}}}%
\end{center}%
\caption{Elementary knot theory, part 2}%
\end{figure}
\end{itemize}

\section{The Duflo isomorphism} \lbl{sec:intro-Duflo} The first of
these equations, ``$\One+\One=\Two$'', is related with the
{\em Duflo isomorphism}.  The Duflo isomorphism is an
algebra isomorphism between the invariant part of the symmetric
algebra and the center of the universal enveloping algebra for any Lie
algebra $\frakg$.  (The Poincaré-Birkhoff-Witt theorem gives a vector
space isomorphism.)  This isomorphism was first described for
semi-simple Lie algebras by Harish-Chandra,
and for all Lie algebras by Duflo~\cite{Duflo:Operateurs}.

Let us review briefly the Duflo isomorphism.
Every Lie algebra $\frakg$ has two associated algebras: the universal
enveloping algebra $U(\frakg)$, generated by $\frakg$ with relations
$xy - yx = [x, y]$, and the symmetric algebra $S(\frakg)$, generated
by $\frakg$ with relations $xy - yx = 0$.  There is a natural map
between the two,
\[
\chi: S(\frakg) \longrightarrow U(\frakg),
\]
given by taking a monomial $x_1\dots x_n$ in $S(\frakg)$ and averaging
over the product (in $U(\frakg)$) of the $x_i$ in all possible
orders. By the Poincaré-Birkhoff-Witt (PBW) theorem, $\chi$ is an
isomorphism of vector spaces and $\frakg$-modules.  Since $S(\frakg)$
is abelian and $U(\frakg)$ is not, $\chi$ is clearly not an algebra
isomorphism.  Even restricting to the invariant subspaces on both
sides,
\[
\chi: S(\frakg)^\frakg \longrightarrow U(\frakg)^\frakg\cong Z(\frakg),
\]
$\chi$ is still not an isomorphism of algebras.  The following theorem
gives a modification that is an algebra isomorphism.

\begin{theorem}[Duflo~\cite{Duflo:Operateurs}]\lbl{thm:duflo} For any
finite-dimensional Lie algebra $\frakg$, the map
\[
\Upsilon: S(\frakg)^\frakg \longrightarrow Z(U(\frakg))
\]
is an algebra isomorphism, where
\begin{equation*}
\begin{split}
\Upsilon &= \chi \circ \partial_{j^{\frac12}} \\
j^{\frac12}(x) &= 
 \det\nolimits^{\frac12}\left(
  \frac{\sinh(\frac{1}{2} \ad x)}
       {\frac{1}{2} \ad x}\right)
\end{split}
\end{equation*}
\end{theorem}

The notation of $\partial_{j^{\frac12}}$ means to consider
$j^{\frac12}(x)$ as a power series on $\frakg$ and plug in the
(commuting) vector fields $\partial/\partial x^*$ on $\frakg^*$.
(Note that for $x^*\in\frakg^*$, $\partial/\partial x^*$ transforms
like an element of $\frakg$).  The result is an infinite-order
differential operator on $\frakg^*$, which we can then apply to a
polynomial on $\frakg^*$ ($\equiv$ an element of $S(\frakg)$).
$j^\frac{1}{2}(x)$ is an important function in the theory of Lie
algebras.  Its square, $j(x)$, is the Jacobian of the exponential
mapping from $\frakg$ to the Lie group $G$.

\section{Wheels} The bridges between the knot theory of
Section~\ref{sec:ElemKnotTheory} and the seemingly quite disparate Lie
algebra theory of Section~\ref{sec:intro-Duflo} are a certain spaces
of uni-trivalent diagrams modulo local relations.  (``Uni-trivalent''
means that the vertices have valence 1 or 3; the 1-valent vertices are
called the ``legs'' of the diagram.)  On one hand, diagrams give
elements of $U(\frakg)$ or $S(\frakg)$ for every metrized Lie algebra
$\frakg$ in a uniform way, as we will see in
Chapter~\ref{chap:liealg}; on the other hand, they occur naturally in
the study of finite type invariants of knots, as we will see in
Chapter~\ref{chap:vassiliev}.  Like the associative algebras
associated to Lie algebras, these diagrams appear in two different
varietes: $\cA$, in which the legs have a linear order, as in
Figure~\ref{fig:A-examp}, and $\cB$, in which the legs are unordered,
as in Figure~\ref{fig:B-examp}.  As for Lie algebras, they each have a
natural algebra structure (concatenation and disjoint union,
respectively); and there is an isomorphism $\chi:\cB\rightarrow
\cA$ between the two (averaging over all possible orders of the legs).

\begin{figure}[tbp]
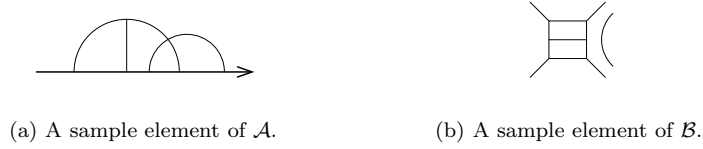
%
\begin{center}%
  \subfigure[A sample element of $\cA$.]%
   {\makebox[.3\linewidth]{\figscaled{A-examp}{0.3}\lbl{fig:A-examp}}}%
\qquad\subfigure[A sample element of $\cB$.]%
   {\makebox[.3\linewidth]{\figscaled{B-examp}{0.25}\lbl{fig:B-examp}}}%
\end{center}%
\caption{Example web diagrams}%
\end{figure}
There is one element of the algebra $\cB$ that will be particularly
important for us: the ``wheels'' element of the title.  It is the
diagrammatic analogue of the function $j^\frac12$ above.
\begin{equation} \lbl{eq:OmegaDef}
  \Omega=\exp \sum_{n=1}^\infty b_{2n}\omega_{2n} \in \cB.
\end{equation}
where
\begin{itemize}
\item The `modified Bernoulli numbers' $b_{2n}$ are defined by the power
series expansion
\begin{equation} \lbl{eq:MBNDefinition}
  \sum_{n=0}^\infty b_{2n}x^{2n} = \frac{1}{2}\log\frac{\sinh x/2}{x/2}.
\end{equation}
These numbers are related to the usual Bernoulli numbers $B_{2n} =
4n\cdot(2n)!\cdot b_{2n}$ and to the values of the Riemann
$\zeta$-function on the even integers. The first three modified
Bernoulli numbers are $b_2=1/48$, $b_4=-1/5760$, and $b_6=1/362880$.

\item The `$2n$-wheel' $\omega_{2n}$ is the degree $2n$ web diagram
made of a $2n$-gon with $2n$ legs:
\begin{equation}\label{eq:omega-definition}
  \omega_2=\mathcenter{\figscaled{2wheel}{0.5}},\quad
  \omega_4=\mathcenter{\figscaled{4wheel}{0.5}},\quad
  \omega_6=\mathcenter{\figscaled{6wheel}{0.5}},\quad\ldots.
\end{equation}
\end{itemize}

Our main theorem is written in terms of $\partial_\Omega$, the
operation of applying $\Omega$ as a differential operator, which takes
a diagram $D$ and
attaches some of its to $\Omega$.  (See
Section~\ref{sec:diagram-differential} for the precise definition.)

\begin{theorem}[Wheeling; joint with D.~Bar-Natan]\lbl{thm:wheeling}
The map $\Upsilon = \chi \circ \partial_\Omega: \cB \rightarrow \cA$
is an algebra isomorphism.  \end{theorem}

(All the notation above, including the definitions of $\cA$ and $\cB$,
is explained in Chapter~\ref{chap:liealg}.  The proof of the theorem
is in Chapter~\ref{chap:wheeling}.)

Although Theorem~\ref{thm:wheeling} was motivated by Lie algebra
considerations when it was first
conjectured~\cite{BGRT:WheelsWheeling, Deligne:letter}, the proof we
will give, based on the equation ``$1+1=2$'' from
Section~\ref{sec:ElemKnotTheory}, is entirely independent of Lie
algebras and is natural from the point of view of knot theory.  In
particular, we obtain a new proof of Theorem~\ref{thm:duflo} for
metrized Lie algebras, with some advantages over the original proofs
by Harish-Chandra, Duflo, and Cartan: our proof does not require any
detailed analysis of Lie algebras, and so works in other contexts in
which there is a Jacobi relation.  For instance, our proof works for
super Lie algebras with no modification.

Theorem~\ref{thm:wheeling} has already seen several applications.  In
Chapter~\ref{chap:wheels} we will use it to compute the Kontsevich
integral of the unknot, using our second elementary knot theory
identity ``$n\cdot 0 = 0$''.

\begin{theorem}[Wheels; joint with T.~Le]\label{thm:wheels}
The Kontsevich integral of the unknot is
\[
Z(\bigcircle) = \Omega \in \cB.
\]
\end{theorem}

We also compute the Kontsevich integral of the
Hopf link $\Hopf$; which is intimately related to the
map $\Upsilon$ above.

\section{Related work} Theorem~\ref{thm:wheeling} was first
conjectured by Deligne~\cite{Deligne:letter} and Bar-Natan,
Garoufalidis, Rozansky, and Thurston~\cite{BGRT:WheelsWheeling}, who
also conjectured Theorem~\ref{thm:wheels}.

Further computations for a sizeable class of knots, links, and
3-manifolds (including all torus knots and Seifert-fiber homology
spheres) have been done by Bar-Natan and
Lawrence~\cite{Bar-NatanLawrence:RationalSurgery}.  Hitchen and
Sawon~\cite{HitchinSawon:CurvatureHyperKahler} have used
Theorem~\ref{thm:wheeling} to prove an identity expressing the $L^2$
norm of the curvature tensor of a hyperk\"ahler manifold in terms of
Pontryagin classes.  And in a future paper~\cite{Thurston:Torus} I
will show how to write simple formulas for the action of $\sltwo(\ZZ)$
on the vector space associated to a torus in the perturbative TQFT of
Murakami and Ohtsuki~\cite{MO:TQFTUniversal}.

These applications suggest that Theorem~\ref{thm:wheeling} is a
fundamental fact.  Another sign of fundamental facts is that they have
many proofs.  Besides the earlier work of Harish-Chandra, Duflo, and
Cartan mentioned above, there are two other recent proofs of
Theorem~\ref{thm:wheeling}. One is due to
Kontsevich~\cite[Section~8]{Kontsevich:DeformationQuantization}, as
expanded by~\cite{ADS:KontsevichQuantization}.  Kontsevich's proof is
already at a diagrammatic level, similar to the one in this thesis,
although it is more general: it works for all Lie algebras, not just
metrized ones.  His proof again uses a transcendental integral,
similar in spirit to the one in Chapter~\ref{chap:kontsevich}.
Another proof is due to Alekseev and
Meinrenken~\cite{AlekseevMeinrenken:NonCommutativeWeil}.  The Alekseev
and Meinrenken paper is not written in diagrammatic language, but
seems to extend to the diagrammatic context without problems.  Their
proof does not involve transcendtal integrals: the only integral in
their proof is in the proof of the Poincaré lemma (the homology of
$\RR^n$ is trivial in dimension $>\,0$).

\section{Plan of the thesis} The first few chapters are standard
introductory material, included here to make this thesis as
self-contained as possible.  In Chapter~\ref{chap:liealg}, we review
the diagrammatics of Lie algebras and define the space $\cA$ and
$\cB$.  In Chapter~\ref{chap:vassiliev} we review the definition of
Vassiliev invariants, again arriving at the same spaces $\cA$ and
$\cB$; in addition, a few more spaces are defined in
Section~\ref{sec:diagrams2}.  The exposition is related to the
standard exposition~\cite{BarNatan:OnVassiliev} and the work of
Goussarov~\cite{Goussarov:KnottedGraphs, Goussarov:3Manifolds} and
Habiro~\cite{Habiro:Claspers}, but is different from both.
Chapter~\ref{chap:kontsevich} is a review of the Kontsevich integral
and its main properties, notably the behaviour under connect sum and
cabling.

After this introductory material, we give the proof of
Theorem~\ref{thm:wheeling} in Chapter~\ref{chap:wheeling}, modulo
fixing the coefficients $b_{2n}$, which is done in
Appendix~\ref{app:coefficients}.  In Chapter~\ref{chap:wheels}, we
prove Theorem~\ref{thm:wheels}, as well as similar results for the
Hopf link.

\chapter{Diagrammatics of Lie algebras}
\lbl{chap:liealg}
Throughout this thesis, we will work heavily with a graphical notation
for tensors, in particular inside a Lie algebra; in fact, this becomes
more than just a notation for us.  For the benefit of readers who may be
unfamiliar with this notation, here is a quick refresher course.

\section{Elementary tensors}
\begin{figure}[hbtp]
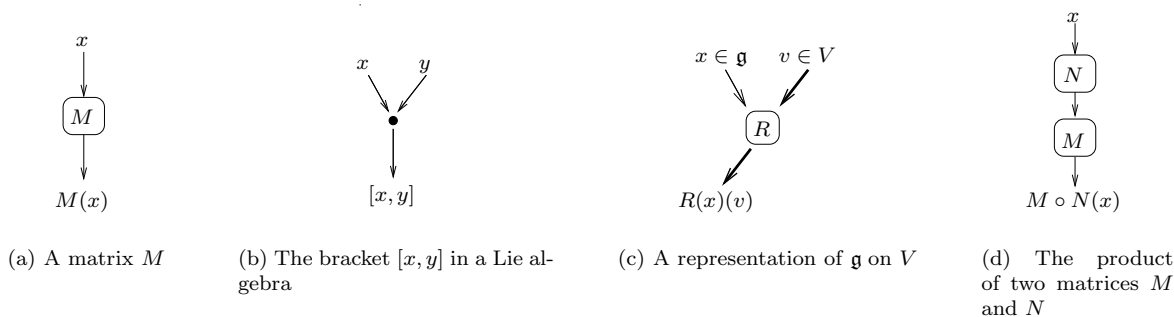
%
\begin{center}%
  \subfigure[A matrix $M$]%
   {\makebox[.19\linewidth]{\figscaled{matrix-simple}{0.55}}
    \label{fig:elemmatrix}}%
\hfil\subfigure[The bracket $\lbracket x,y\rbracket$ in a Lie algebra]%
   {\makebox[.29\linewidth]{\figscaled{bracket}{0.55}}
    \label{fig:elembracket}}%
\hfil\subfigure[A representation of $\frakg$ on $V$]%
   {\makebox[.29\linewidth]{\figscaled{representation}{0.55}}
    \label{fig:elemrep}}%
\hfil\subfigure[The product of two matrices $M$ and $N$]%
   {\makebox[.19\linewidth]{\figscaled{matrix-mult}{0.55}}
    \label{fig:matrixmult}}%
\end{center}%
\caption{Some elementary tensors}\lbl{fig:elemtensors}%
\end{figure}

In general, a tensor with $n$ indices will be represented by a graph
with $n$ legs (or free ends).  The indices of a tensor can belong to
various vector spaces or their duals.  Correspondingly, the legs of
the graph should be colored to indicate the vector space and directed
to distinguish between a vector space and its dual.  See
Figure~\ref{fig:elemtensors} for some examples.
Figure~\ref{fig:elemmatrix} shows a generic matrix $M \in \Hom(V, V)
\isom V^*\otimes V$.  By convention, data flows in the direction of
the arrows, so the incoming edge is the $V^*$ factor and the outgoing
arrow is the $V$ factor.  Figure~\ref{fig:elembracket} represents the
bracket in a Lie algebra $[\cdot,\cdot] \in \frakg^* \otimes \frakg^*
\otimes \frakg$.  Because these vertices are so ubiquituous, these
vertices will not be decorated.  A representation $R$ on a vector
space $V$ of a Lie algebra is a linear map $R \in \Hom(\frakg \otimes
V, V) \equiv \frakg^* \otimes V^*n
\otimes V$.  This is depicted in Figure~\ref{fig:elemrep}.

\section{Composing Tensors}
The use of tensors is that you can compose them in many different
ways.  For instance, two matrices $M, N \in Hom(V, V)$ can be
multiplied:
\[
M \circ N = \sum_\beta M^\alpha_\beta N^\beta_\gamma
\]
which we can represent graphically as in Figure~\ref{fig:matrixmult}.
More generally, if a tensor has two indices with values in vector
spaces that are dual to each other, they can be contracted;
graphically, an incoming and outgoing leg of the same color can be
connected.

As an important example, the
Jacobi relation in a Lie algebra,
\[
[x,[y,z]] + [z,[x,y]] + [y,[z,x]] = 0, \qquad x,y,z \in \frakg
\]
can be expressed graphically as 
\begin{center}
\figscaled{jacobi-directed}{0.4}
\end{center}
In the same way, the equation that a representation of a Lie algebra
be a representation can be written graphically.
\begin{align*}
R([x,y]) &= R(x)R(y) - R(y)R(x) \\
\figcent{rep-bracket}{0.4} &= \figcent{rep-xy}{0.4} - \figcent{rep-yx}{0.4}
\end{align*}

As another important case, suppose we have a metrized Lie algebra:
that is, we have an invariant nondegenerate bilinear form
on $\frakg$.  This gives us a tensor $(\cdot\comma\cdot) \in \frakg^*
\otimes \frakg^*$ and its inverse in $\frakg \otimes \frakg$.
Graphically, this gives us tensors
\[
\figcent{metric}{0.5} \quad \text{and}\quad \figcent{metric-inv}{0.5}
\]
With the aid of these two tensors, we can glue two edges on which the
arrows don't match.  Furthermore, because they are inverses of each
other, it doesn't matter if we stick in extra pairings.  We will take
this as license to drop all arrows on Lie algebra legs when working
with metrized Lie algebras.

There is one point to be careful about: when we drop the decorations
on a Lie bracket, we are assuming some symmetry of the bracket, since
it is impossible to tell the bracket from a rotated version of itself.
Fortunately, the Lie bracket is cyclicaly invariant when you identify
$\frakg$ with $\frakg^*$ using $(\cdot\comma\cdot)$ by invariance of
the metric.

\section{The function $j^{\frac12}$}\lbl{sec:function-j} To gain
practice with the graphical notation for tensors, let us find the
graphical version of the function $j^{\frac12}$ on $\frakg$.  Let us
recall the definition of $j^{\frac12}$:
\[
j^{\frac12}(x) = \det^{\frac12} \frac{\sinh \frac {\ad x} 2}{\frac {\ad x} 2}.
\]
We will work through this definition step by step.

\begin{figure}[tbp]
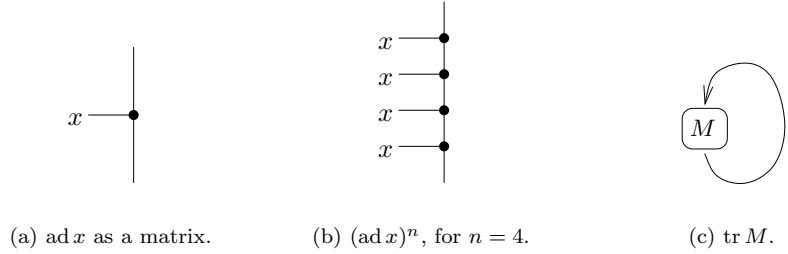

\begin{center}
\subfigure[$\ad x$ as a matrix.]%
 {\makebox[.25\linewidth]{\figscaled{adx}{0.6}\lbl{fig:adx}}}%
\subfigure[$(\ad x)^n$, for $n=4$.]%
 {\makebox[.25\linewidth]{\figscaled{adx-4}{0.6}\lbl{fig:adxn}}}%
\subfigure[$\tr M$.]
 {\makebox[.25\linewidth]{\figscaled{matrix-trace}{0.6}\lbl{fig:matrix-trace}}}%
\end{center}
\caption{Building blocks of wheels}
\end{figure}
\begin{itemize}
\item $(\ad x)^n$.
Here we take $\ad x$, considered as a matrix acting on $\frakg$, and
raise it to some power.  Since $(\ad x)(y) = [x,y]$, the graphical
representation of $\ad x$ is as in Figure~\ref{fig:adx}.
Note that the end labelled $x$ is not considered as a leg in the tensor sense;
instead, we plug in the (fixed) element $x \in \frakg$ and get a fixed
matrix.  But see Section~\ref{sec:symalg} for another point of view.

To raise $\ad x$ to a power, just take several copies and string them
together as in Figure~\ref{fig:adxn}.

\item Determinants and traces.
The determinant of a matrix is a non-linear function of
the matrix that does not fit well in our graphical notation.
Fortunately, we can get rid of determinants using the equality
\[
\det \exp M = \exp \tr M.
\]
The trace of a matrix is easy to understand graphically: take a matrix
$M^a_b \in V^* \otimes V$ and sum over the diagonal, $a=b$; in other
words, contract $V$ with $V^*$.  Graphically, we just connect the
input of $M$ with its output as in Figure~\ref{fig:matrix-trace}. 
\end{itemize}

We can now write a complete formula for $j^{\frac12}$ in terms of
graphs:
\ssp
\begin{equation}
\begin{split}
j^{\frac12}(x) &= \det^{\frac12}\frac{\sinh \frac {\ad x} 2}{\frac {\ad x} 2}\\
	&= \exp(\frac12\tr(\log \frac{\sinh \frac {\ad x} 2}{\frac {\ad x} 2}))\\
	&= \exp(\frac12\tr(\sum_{n=0}^\infty b_{2n} (\ad x)^{2n}))\\
	&= \exp(\sum_{n=0}^\infty b_{2n} \omega_{2n}(x)) \\
	&= \exp\left(\frac{1}{48}\figcent{2wheelx}{0.3}
		   - \frac{1}{5760}\figcent{4wheelx}{0.3}
                   + \frac{1}{362880}\figcent{6wheelx}{0.3} - \cdots\right)
\end{split}
\end{equation}
\dsp
The $b_{2n}$ were defined in Equation~\ref{eq:MBNDefinition}.  The
$\omega_{2n}$ are as in Equation~\ref{eq:omega-definition}, with $x$
placed on the legs.

\section{The symmetric and universal enveloping algebras}
\lbl{sec:symalg}
All the tensors we have drawn so far have had a fixed number of legs.
We will also be interested in representing diagrammatically the
symmetric algebra of a Lie algebra or its universal enveloping algebra
so that we can, e.g., write $j^{1/2}(x)$ as an element of $S(\frakg)$
rather than just a function on $\frakg$.

Recall that
\[
S(\frakg) = \bigoplus_{n \ge 0} \Sym^n(\frakg)
\]
where $\Sym^n(\frakg)$ is the $n$'th symmetric power of $\frakg$,
which we take to be the quotient of $\frakg^{\otimes n}$ by the
symmetric group $S_n$.  Diagrammatically, an element of
$\Sym^n(\frakg)$ is easy to represent; just take a diagram with $n$
legs (representing an element of $\frakg^{\otimes n}$) and take the
quotient by the symmetric group, i.e., forget about the labels on the
legs.  An element of $S(\frakg)$ is a diagram with any number of
unlabelled legs.

For the universal enveloping algebra $U(\frakg)$, recall the
definition: $U(\frakg)$ is the associative algebra generated by
$\frakg$, modulo the relation $xy - yx = [x,y]$ for $x,y \in\frakg$.
In particular, it is a quotient space of the free associative algebra
on $\frakg$, which is $\bigoplus_n \frakg^{\otimes n}$.  These are
diagrams with $n$ legs, with certain relations.  We will remember
these relations by placing the legs on a thick line, as in
Figure~\ref{fig:A-examp}.  (Explicitly, we evaluate the interior part of
the diagram to get an element in $\frakg^{\otimes n}$; we then
multiply the $n$ elements of $\frakg$ to get an element of
$U(\frakg)$.  The quotient means that the relation 
\[
\mlarge\SGraph = \mlarge\TGraph- \mlarge\UGraph
\]
is satisfied.

\section{Spaces of diagrams}\lbl{sec:diagrams1}

Until now, the diagrams we have drawn have only been a particular
notation for tensors in Lie algebras.  In fact, the notation is more
than a notation: they will be the fundamental objects of study.  We
will now define a space $\cB$ that is informally the space of all
tensors in $S(\frakg)$ that can be constructed from the Lie
bracket, modulo the relations that can be deduced solely from the
Jacobi identity.
n
\begin{definition} An {\em open Jacobi diagram} (variously called a
Chinese Character, uni-trivalent graph, or web diagram) is
vertex-oriented uni-trivalent graph, i.e., a graph with univalent and
trivalent vertices together with a cyclic ordering of the edges
incident to the trivalent vertices.  Self-loops and multiple edges
are allowed.  The univalent vertices are called {\em legs}.  In planar
pictures, the orientation on the edges incident on a vertex is the
clockwise orientation.  Some examples are shown in
Figure~\ref{fig:B-examp}.
\end{definition}

\begin{definition} $\cB^f$ is the vector space spanned by Jacobi diagrams
modulo the IHX relation {$\IGraph = \HGraph - \XGraph$} and the
anti-symmetry relation $\YGraph + \TwistedYGraph = 0$, which can be
applied anywhere within a diagram.  \end{definition}

\begin{definition} The {\em degree} of a diagram in $\cB^f$ is half
the number of vertices (trivalent and univalent).  $\cB$ is the
completion of $\cB^f$ with respect to the grading by degree.
\end{definition}

By the discussion above, any element of $\cB^f$ gives an element in
$S(\frakg)$ for any metrized Lie algebra $\frakg$.  In some sense, you
can think of $\cB$ as being related to a ``universal (metrized) Lie
algebra'', incorporating information about all Lie algebras at once.
But $\cB$ is both bigger and smaller than that: the map to the
product of $S(\frakg)$ for all metrized Lie algebras is neither
injective nor surjective.
\begin{itemize}
\item There are elements of $\cB$ that are non-zero but
      become zero when evaluated in any metrized Lie algebra.
      See Vogel~\cite{Vogel:Structures} for details.
\item Not all elements of $S(\frakg)$ are in the image
      of the map from $\cB$ to $S(\frakg)$, as shown by
      the following lemma:
\begin{lemma} Every element in the image of $\cB$ is invariant under
the action of any Lie group $G$ (not necessarily connected) whose Lie
algebra is $\frakg$.
\end{lemma}
\begin{proof}
Every element in the image of $\cB$ is made by contracting copies of
the structure constants of the Lie algebra, which are invariant under
$G$.
\end{proof}
\end{itemize}

There is a natural algebra structure on $\cB$, the disjoint union
$\cup$ of diagrams, which corresponds to the algebra structure on
$S(\frakg)$.

We have similarly a diagrammatic analogue of $U(\frakg)$:

\begin{definition} A {\em based Jacobi diagram} (also called chord
diagram or Chinese Character diagram) is a Jacobi diagram with a
total ordering on its legs.  They are conventionally represented as
in~\ref{fig:A-examp}.
\end{definition}

\begin{definition} $\cA^f$ is the vector space of
based Jacobi diagrams modulo the Jacobi and antisymmetry relations as in
$\cB^f$, plus the STU relation $\SGraph = \TGraph - \UGraph$.
\end{definition}

\begin{definition}
The {\em degree} of a diagram in $\cA^f$ is half the number of
vertices.  $\cA$ is the completion of $\cA^f$ with respect to the
grading by degree.
\end{definition}

There is likewise a natural algebra structure on $\cA$: take two based
Jacobi diagrams $D_1$, $D_2$ and place the legs of $D_1$ before the legs
of $D_2$ in the total ordering on legs.  For reasons that will become
clear later, we call this the {\em connected sum} and denote it $D_1 \connect
D_2$.

These two products live on isomorphic spaces so may be confused.  We
usually write out the product in cases of ambiguity.  If an explicit
symbol for the product is omitted, the product is the disjoint union
product $\udot$ unless otherwise specified.

\section{Diagrammatic Differential Operators}
\lbl{sec:diagram-differential}
Our main goal in this thesis is to find a diagrammatic analogue of
the Duflo isomorphism of Theorem~\ref{thm:duflo}.  In that formula the
differential operator $\partial_{j^{1/2}}$
plays a prominent role.  We saw above how to write $j^{\frac12}$ as an
element of $\cB$; but what does it mean to apply it as a differential
operator?

\begin{definition}\label{def:diagram-differential} For a diagram $C
\in \cB$ without struts (components like $\strutn{}{}$), the operation
of {\em applying C as a differential operator}, denoted $\partial_C: \cB
\rightarrow \cB$, is defined to be
\ssp
\[ \partial_C(D)=\begin{cases}
    0 & \parbox{1.7in}{if $C$ has more legs than $D$,} \\[3pt]
    \parbox{2.2in}{
      the sum of all ways of gluing all the legs of $C$ to some (or all)
      legs of $D$
    }\quad & \text{otherwise.}
  \end{cases}
\]
\dsp
\end{definition}

For example,
\[ \partial_{\omega_4}(\omega_2)=0; \qquad
   \partial_{\omega_2}(\omega_4)=
     8\mathcenter{\figscaled{SideGluing}{0.5}}
	+4\mathcenter{\figscaled{DiagonalGluing}{0.5}}.
\]
If $C$ has $k$ legs and total degree $m$, then $\partial_{C}$ is an operator
of degree $m-k$. By linear extension, we find that every $C\in\cB$
without struts defines an operator $\partial_{C}:\cB\to\cB$.  (We
restrict to diagrams without struts to avoid circles arising from the
pairing of two struts and to guarantee convergence: gluing with a
strut lowers the degree of a diagram, and so the pairing would not
extend from $\cB^f$ to $\cB$.)

We leave it to the reader to verify that this operation $\partial_ C$ is
the correct diagrammatic analogue of applying $C$ as a differential
operator.  (More precisely, these are constant coefficient
differential operators; for instance, they all commute.)  There are
some good signs that we have the right definition:
\begin{itemize}
\item A diagram $C$ with $k$ legs reduces the number of legs by $k$,
      corresponding to a differential operator of order $k$.
\item If $k = 1$ ($C$ has only one leg), we have a Leibniz rule like
      that for linear differential operators:
\[
\partial_C(D_1 \cup D_2) = \partial_C(D_1) \cup D_2 + D_1 \cup\partial_C(D_2).
\]
      (Actually, all diagrams with only one leg in $\cB$ are 0, so
      we have to extend our space of diagrams slightly for this
      equation to be non-empty.  Adding some extra vertices of
      valence 1 satisfying no relations is sufficient.)
\item Multiplication on the differential operator side is the same
      thing as composition:
\[
\partial_{C_1 \cup C_2} = \partial_{C_1} \circ \partial_{C_2}.
\]
\end{itemize}
With a little more work, these properties can give a proof that this
is the correct diagrammatic analogue of differential operators:
extend the space of diagrams as suggested in the second
point and add an operation of gluing two of these univalent ends.  We
will not give the details here.

\section{The map $\Upsilon$}\label{sec:upsilon}
We are now in a position to find the
diagrammatic analogue of the Duflo map $\Upsilon$ we introduced for
Lie algebras.  Recall that $\Upsilon$ is the composition of two maps,
the infinite-order differential operator $\partial_{j^{1/2}}$ and he
Poincaré-Birkhoff-Witt isomorphism $\chi: S(\frakg) \rightarrow U(g)$.
Both of these have natural diagrammatic analogues.

From the computations in Section~\ref{sec:function-j}, the
      diagrammatic analogue of $j^{\frac12}$ is the ``wheels'' element
\[
\Omega = \exp(\sum_{m=0}^\infty b_{2m} \omega_{2m}) \in \cB.
\]
      $\partial_\Omega$ is the action of
      applying this element as a differential operator, as in the
      previous section; this is called the ``wheeling'' map.

The Poincaré-Birkhoff-Witt isomorphism is the isomorphism from
$S(\frakg)$ to $U(\frakg)$ given by taking a monomial $x_1\cdots
x_n$ in $S(\frakg)$ to the average over the product in
$U(\frakg)$ of the $x_i$ in all possible orders.  We can define
a similar map diagrammatically.
\begin{definition}\label{def:chi}
The map $\chi: \cB \rightarrow \cA$ on a Jacobi
diagram $D \in \cB$ with $n$ (unordered) legs is the average (not sum) in
$\cA$ over all possible orders on the legs.
\end{definition}
As for Lie algebras, this is a vector space isomorphism.  There
are several different proofs of this fact.  See
Bar-Natan~\cite[Section~5.2, Theorem~8]{BarNatan:OnVassiliev}
for one.

The natural analogue of Theorem~\ref{thm:duflo} (the Duflo
isomorphism) is that the composition of these two maps be an algebra
isomorphism; this is exactly Theorem~\ref{thm:wheeling}:

\setcounter{temp}{\value{theorem}}
\setcounter{theorem}{1} \begin{theorem}[Wheeling; joint with
D.~Bar-Natan] The map $\Upsilon = \chi \circ \partial_\Omega: \cB
\rightarrow \cA$ is an algebra isomorphism.  \end{theorem}
\setcounter{theorem}{\value{temp}}

Note that $\partial_\Omega$ decreases or keeps fixed the number of legs of
a diagram, and since $\Omega$ starts with 1, $\hat\Omega$ is
lower-triangular with respect to the number of legs and so a vector
space isomorphism.  Since $\chi$ is also a vector space isomorphism,
$\Upsilon$ is automatically bijective.  The content of the theorem is that
$\Upsilon$ is an algebra map:
\begin{equation}\lbl{eq:wheeling}
\Upsilon(D_1 \cup D_2) = \Upsilon(D_1) \connect \Upsilon(D_2).
\end{equation}

\chapter{Vassiliev Invariants}\label{chap:vassiliev}

The same spaces of diagrams that we saw appearing from Lie algebras in
the last chapter (the spaces $\cA$ and $\cB$) also arise from
3-dimensional topology, via {\em Vassiliev} or {\em finite type}
invariants.  In this chapter we review this theory.  Our presentation
is slightly non-standard.  It is related to the standard presentation
of Bar-Natan~\cite{BarNatan:OnVassiliev} as well as the claspers or
clovers of Goussarov~\cite{Goussarov:KnottedGraphs,
Goussarov:3Manifolds} and Habiro~\cite{Habiro:Claspers}.
But the differences are purely expositional; a reader familiar with
the theory should feel free to skip this chapter after reviewing the
definitions in Section~\ref{sec:diagrams2}.

All knots, tangles, etc. in this thesis are framed.  The space $\cK$
is the vector space (over $\CC$ for convenience, though any field of
characteristic 0 would work for this chapter) spanned by framed knots,
i.e., smooth embeddings $S^1 \hookrightarrow \RR^3$ with a normal
vector field at each point of $S^1$, considered up to isotopy.  More
generally, suppose $X$ is a compact 1-manifold, possibly with
boundary.  Then $\cK(X)$ is the vector space spanned by framed tangles
with skeleton $X$, i.e., framed smooth proper embeddings
$X\hookrightarrow B^3$ with some fixed (but usually unspecified)
behaviour at the boundary, considered up to isotopy relative to the
boundary.  In drawings, framings will be the blackboard framing: the
framing perpendicular to the plane of the paper.

\section{Definition of Vassiliev invariants} The idea of finite type
invariants is to find invariants of tangles analagous to polynomial
functions on vector space.  A polynomial $p$ of degree $n$ on a vector
space $V$ can be defined to be a function that vanishes when you take the
alternating sum of its values on the vertices of a parallelepiped:
\[
\text{$p$ of degree $n$} \Leftrightarrow
  \sum_{T \subset S} (-1)^{|T|} p(x_0 + \sum_{v \in T} v) = 0,
	\qquad x_0 \in V, S \subset V, \lvert S\rvert = n+1.
\]
Similarly, a finite type invariant of degree $n$ is a
knot invariant that vanishes when you take the alternating sum on the
vertices of a ``cube of knots'', determined by picking a knot
projection and flipping some subset of the crossing.  

\begin{definition}[First version]\label{def:finite-type-first} A
knot invariant  $f: \cK\rightarrow \CC$%
\footnote{More generally, we can allow values in an arbitrary abelian group}
is {\em finite type} of degree $n$ if, for every knot projection $K$
with and subset $S$ of the crossings and arcs of $K$, $|S| = n+1$,
\[
\sum_{T \subset S} (-1)^{|T|} f(K_T) = 0,
\]
where $K_T$ is the knot obtained from $K$ by flipping the crossings in
$T$ and adding a positive full twist to the framing along the arcs in
$T$.
\end{definition}

We need to allow the framing changes as a special case because
crossing changes can only change the framing by two full twists.  This
definition, like all others in this section, extends
without change to $X$-tangles for an arbitrary 1-manifold $X$.

Vassiliev invariants of degree 0 are constansts.  In degree 1, we
have the {\em linking number}, which is an
invariant of 2-component links $X = \bigcircle_a\bigcircle_b$.
Intuitively, this measures the number of times the component $a$ winds
around the component $b$.  There are several ways to define the
linking number precisely.  One is to count, with signs, the number of
times the component $a$ crosses over the component $b$ in a given
projection.  More intrinsically, there is a map from $S^1\times S^1$
to $S^2$ which takes a point in $S^1\times S^1$ to the normalized
vector pointing from the point on $a$ to the point on $b$.  The
linking number is the degree of this map.  When you change a crossing
in a projection, the linking number changes by 0 or $\pm 1$, depending
only on the intrinsic topology of the crossing (whether the crossing
was a crossing of a component with itself or between
two different components); when you change a second crossing, this
difference doesn't change, so the linking number is finite type of
degree 1.

As a related example, every knot has a canonical framing given by the
boundary of a Seifert surface; the difference between the given
framing of a knot and its canonical framing is again finite type of
degree 1.

Definition~\ref{def:finite-type-first} is slightly unsatisfactory
since it depends on picking a knot projection, which squashes 3
dimensions down to 2 in a slightly uncomfortable fit.  We can fix this
by ``delocalising'' the crossing change.

\begin{definition}\lbl{def:circus} An {\em
$n$-circus} for a knot $K$ is an embedding in $S^3\setminus K$ of $n$
disjoint {\em single lassos} resp. {\em double lassos}
\begin{center}
\figcent{single-lasso}{0.3}\quad\text{resp.}\quad\figcent{double-lasso}{0.3},
\end{center}
that are embedded trivially in $S^3$ in the absence of $K$.
\end{definition}

\ssp
\begin{figure}[tbp]
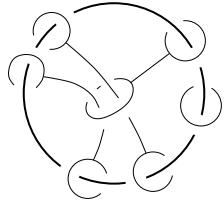

\begin{center}
\figscaled{nCircus-examp}{0.40}
\end{center}
\caption{A sample 4-circus.}\lbl{fig:circus-ex}
\end{figure}
\dsp

\begin{definition}\label{def:lasso-action} The {\em action}
$\Act_C$ on an $n$-circus $C$ is the move
\begin{equation}\label{eq:dbl-lasso-action}
\Act\left(\figcent{simple-lasso}{0.2}\right) =
	\figcent{simple-crossing}{0.2}
\end{equation}
on each of its component double lassos and the move
\begin{equation}
\Act\left(\figcent{simple-single-lasso}{0.3}\right) =
	\figcent{simple-twist}{0.3}
\end{equation}
on each of its component single lassos.
\end{definition}

The geometric picture for the double lasso action is that the double
lassos is like a kayak paddle.  For those unfamiliar with kayaking, a
kayak paddle has a blade\footnote{The part of a paddle that goes in
the water.} on either end.  The blades at the two ends are at a
$90^\circ$ angle to each other.  For us, the loops run around the
outside of the blades the strands pass perpendicularly through the
loops.  Then you can shrink the handle down to zero and pass the
strands through each other with no ambiguity, since they meet at a
right angle.  Unfortunately, this is hard to draw in a 2-dimensional
picture.  For the drawings in Equation~\ref{eq:dbl-lasso-action} to
really correspond to this geometric picture, we would need to include
an extra (negative) quarter-twist in the handle of the lasso, but this
would clutter future diagrams too much.  The conventions here agree
with those of Habiro~\cite{Habiro:Claspers}, but are opposite those of
Goussarov~\cite{Goussarov:KnottedGraphs,Goussarov:3Manifolds}.

\begin{definition}[Second version]\lbl{def:fintype} A knot invariant
$f$ is {\em finite type} of degree $n$ if, for every $n+1$-circus $S$,
\[
\sum_{T \subset S} (-1)^{|T|} f(\Act_T(K)) = 0.
\]
The formal linear combination of knots $\sum_{T \subset S} (-1)^{|T|}
\Act_T(K)$ is called the {\em resolution} $\delta^{(n+1)}(S)$ of $S$.
\end{definition}

\begin{remark}
A double lasso can be made out of three single lassos:
\[
\Act\left(\figcent{simple-lasso}{0.2}\right) =
	\Act\left(\figcent{three-single-lassos}{0.2}\right)
\]
(where the ``$-1$'' indicates a single lasso to be applied in the
reverse sense), so we could have have omitted double lassos from
Definition~\ref{def:circus} without changing the definition of finite type.
But we prefer not to do this: instead, we will write single
lassos in terms of double lassos.  (See Relation~\ref{rel:sing-lasso}.)
The relationship between single and double lassos is comparable to
the relationship between quadratic and bilinear forms.
\end{remark}

\section{Weight systems}\lbl{sec:weight-systems} The resolution of an
$n+2$-circus is the difference of the resolutions of two $n+1$-%
circuses, so a Vassiliev invariant of degree $n$ is also a Vassiliev
invariant of degree $m$ for $m>n$; i.e., Vassiliev invariants form an
increasing filtration.  The image of a Vassiliev invariant in the
associated graded space is called its {\em weight system}; it is
comparable to the highest degree term of a polynomial or the symbol of
a differential operator.  In this section we will try to identify this
associated graded space.  It is more convenient to work with the dual
picture.

\begin{definition} The $n$'th term of the {\em Vassiliev filtration}
$\cK = \cK_0\superset \cK_1
\superset\cdots\superset\cK_n\superset\cdots$ on the space of knots is
the span of resolutions of $n$-component
circuses in the complement of a knot.
\end{definition}

Let $\Diag(C)$ or the {\em diagram} associated to $C$ be the image of
$\delta^{(n)}(C)$ in $\cK_n/\cK_{n+1}$.

\begin{relation}[Homotopy]\lbl{rel:homotopy} $\Diag(C) \equiv
\Diag(C')$, where $C$ and $C'$ differ by homotopy of the knot and
handles of the lasso (fixing the endpoints) in the complement of the
loops of lassos.
\end{relation}

\begin{proof} Given an $n$-circus $C$, consider the $n+1$-circus
$\tilde C$ with an additional double lasso which encircles two knot
strands and/or handles of double lassos.  The resolution of $\tilde C$
along the new double lasso is equal to the difference between $C$ and
a variant $C'$ in which the two strands have crosses each other.
Since $\tilde C$ has $n+1$ lassos, this difference is zero in the
associated graded space $\cK_n/\cK_{n+1}$ so $\Diag(C) = \Diag(C')$.
Likewise framing changes on a strand can be achieved by adding an
additional single lasso looping it.  These two moves generate
homotopy as in the statement.
\end{proof}

With the Homotopy Relation, we have forgotten much of the topology of
the embedding of an $n$-circus.  The next relation lets us simplify
the loops of double lassos.

\begin{relation}[Splitting]\lbl{rel:splitting} 
Let $C, C_1, \dots, C_k$ be $n$-circuses which are the same except for
one loop of one double lasso.  For that loop, suppose that the loop of
$C$ is homotopically the sum of the loops of $C_1, \dots, C_k$.  Then
$\Diag(C) = \Diag(C_1) + \dots + \Diag(C_k)$.
\end{relation}

\begin{proof}
A multiple crossing can be done step by step:
\[
\delta\left(\mathcenter{\fig{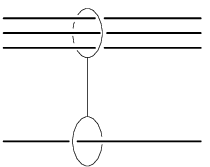}{0.25}}\right) =
	\delta\left(\mathcenter{\fig{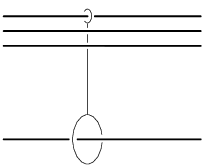}{0.25}}\right)
	+ \delta\left(\mathcenter{\fig{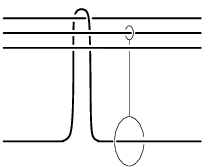}{0.25}}\right)
	+ \delta\left(\mathcenter{\fig{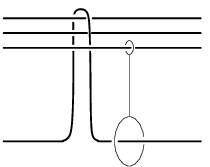}{0.25}}\right).
\]
($\delta$ is the operation of resolving the double lasso.)  Modulo
Relation~\ref{rel:homotopy}, we can forget the extra little hooks.
\end{proof}

\begin{relation}[Single lasso]\lbl{rel:sing-lasso} Let
$C$ be a circus with a single lasso and $C'$ be the same circus with a
the single lasso replaced by a double lasso whose two loops are
parallel to the old single lasso.  Then
\[
2\Diag(C) = \Diag(C').
\]
\end{relation}
\begin{proof}
\[
2\delta\left(\figcent{single-lasso-simp}{0.3}\right)
= \delta\left(\figcent{single-lasso-simp}{0.3}\right) +
	\delta\left(\figcent{twist-lasso}{0.25}\right)
= \figcent{line}{0.25} - \figcent{double-twist}{0.25}
= \delta\left(\figcent{double-lasso-chord}{0.25}\right).
\]
\end{proof}

\begin{figure}[tbp]
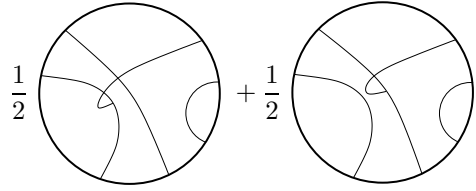
\label{fig:Jacobi-ex}
\[
\frac{1}{2}\figcent{nCircus-diag-1}{0.4} +
\frac{1}{2}\figcent{nCircus-diag-2}{0.4}
\]
\caption{The diagrams corresponding to Figure~\ref{fig:circus-ex}.}
\end{figure}

Let us see what the space of diagrams $\cK_n/\cK_{n+1}$ is, using only
Relations~\ref{rel:homotopy},~\ref{rel:splitting},
and~\ref{rel:sing-lasso}.  By Relation~\ref{rel:sing-lasso}, we can
get turn all single lassos into double lassos.  By
Relation~\ref{rel:splitting}, we can make sure each loop of a double
lasso encloses only a single strand (either a strand of the knot or a
handle of a double lasso); and by Relation~\ref{rel:homotopy} only the
combinatorics of these loops and connections are relevant.  So we can
reduce any diagram $\Diag(C)$ to a sum of encircled, enriched Jacobi
diagrams as in Figure~\ref{fig:Jacobi-ex}.

\ssp
\dsp

\begin{definition} An {\em encircled Jacobi diagram} is a
vertex-oriented trivalent graph with a distinguished cycle, called the
{\em skeleton}.  As a special case, a distinguished circle with no
vertices is allowed.  Vertices on the distinguished circle are called
{\em external}; other vertices are called \emph{internal}.
\end{definition}

\begin{definition}
An \emph{routed Jacobi diagram} is a Jacobi diagram with two of the
vertices incident to each internal vertex distinguished.
\end{definition}

In drawings, the distinguished cycle is indicated by thick lines. Thin
lines are double lassos and dashed lines are either of the two.  For
routed Jacobi diagrams, the two distinguished edges are continuous and
the remaing edge meets them at right angles.

Most routed encircled Jacobi diagrams correspond to
$n$-circuses in a knot complement by embedding the diagram arbitrarily
in $\RR^3$, resolving each vertex like this
\[
\figcent{vertex}{0.5} \rightsquigarrow
	\figcent{vertex-expand}{0.5}.
\]
The distinguished cycle becomes the knot.  In good cases, this
yields an $n$-circus in the knot complement.  We will examine this
situation in Appendix~\ref{app:cyclic}, where we will prove the
following proposition.

\begin{definition}
A Jacobi diagram is {\em boundary connected} if there are no connected
components disjoint from the skeleton (i.e., the distinguished cycle
in the case of encircled Jacobi diagrams).
\end{definition}

\begin{proposition}\label{prop:cyclic} Any boundary connected Jacobi
diagram has a routing which corresponds to an $n$-circus.
Furthermore, the diagrams of all such $n$-circuses are equal.
\end{proposition}

Proposition~\ref{prop:cyclic} gives us license to forget the
routing on Jacobi diagrams so long as we restrict to boundary
connected diagrams.

The relations above have more consequences.  Any diagram containing an
empty loop is 0 by the Splitting relation.  But this empty loop need
not be obviously empty in the associated Jacobi diagrams.  As the loop
is deformed, it can pass over strands or vertices; these give us two
types of relations.  Before we introduce the new relations, let us
define a handy graphical notation.
\begin{definition}\label{def:rbox}
A rounded box like
\[
\figcent{rloop-label}{0.2}
\]
is, by definition, the sum over the $n$ ways of attaching the strand
$s_0$ to one of the strands $s_1$, \dots, $s_n$:
\[
\figcent{rloop}{0.2}
   = \figcent{rloop-1}{0.2}
   + \figcent{rloop-2}{0.2}
   + \ldots
   + \figcent{rloop-n}{0.2}.
\]
\end{definition}

The relations are now
\begin{enumerate}
\item (Antisymmetry) We can push a strand (either a piece of knot or a
      handle) through the loop of an empty lasso:
\begin{equation}\lbl{eq:antisymmetry}
0 = \Diag\bigp{\figcent{AS-1}{0.19}}
  = \Diag\bigp{\figcent{AS-2}{0.19}} + \Diag\bigp{\figcent{AS-3}{0.19}}
\end{equation}
or
\[
0 = \figcent{AS-diag-0}{0.19}
	= \figcent{AS-diag-1}{0.19} +\figcent{AS-diag-2}{0.19}.
\]
In other words, reversing the orientation of a trivalent vertex
negates the diagram.

\item (vertex) We can push a vertex through the loop of a lasso:
\begin{equation}\lbl{eq:vertex}
0 = \Diag\bigp{\figcent{jacobi-1}{0.19}} =
  \Diag\bigp{\figcent{jacobi-2}{0.19}} +
  \Diag\bigp{\figcent{jacobi-3}{0.19}} +
  \Diag\bigp{\figcent{jacobi-4}{0.19}}
\end{equation}
or
\[
\figcent{jacobi-diag-0}{0.19}
= \figcent{jacobi-diag-1}{0.19} + \figcent{jacobi-diag-2}{0.19} + 
	\figcent{jacobi-diag-3}{0.19} = 0.
\]
      As with the antisymmetry relation, this comes in two different
      versions, depending on whether the thin line is the handle of a
      lasso or a strand of the knot.  The relation is also called the
      Jacobi or IHX relation in the first case, and the universal
      enveloping algebra or STU relation in the second case.
\end{enumerate}

Our approximation to the space of diagrams $\cK_n/\cK_{n+1}$ is the
space $\cA_n(\bigcircle)$ of degree $n$ boundary-connected encircled Jacobi
diagrams modulo the antisymmetry and vertex relations.  We have shown
that the map $\Diag:\{\text{$n$-circuses}\}\rightarrow\cK_n/\cK_{n+1}$
factors through a map $\Diag'$ to $\cA_n(\bigcircle)$:
\[
\Diag: \text{$n$-circuses} \overset{\Diag'}{\longrightarrow}
  \cA(\bigcircle) \overset{m}{\longrightarrow} \cK_n/\cK_{n+1}
\]
So $\cA_n(\bigcircle)$ is an upper bound for $\cK_n/\cK_{n+1}$.  The map
$m$ is actually an isomorphism, but this fact seems to be difficult to
prove directly.  It follows from the construction in
Chapter~\ref{chap:kontsevich}.

\section{Spaces of diagrams, redux}\lbl{sec:diagrams2}

As mentioned previously, all constructions in the previous two
sections extend without change to the more general context of tangles
with skeleton $X$, for $X$ any 1-manifold.  Repeating the discussion of
Section~\ref{sec:weight-systems}, we are naturally led to following
definition for the corresponding spaces of diagrams.

\begin{definition} For a compact 1-manifold $X$ (possibly with boundary), a
{\em Jacobi diagram based on $X$} is a vertex-oriented uni-trivalent
graph $\Gamma$ together with an embedding $\phi: X \hookrightarrow
\Gamma$ up to isotopy so that $\phi(\partial X)$ is
exactly the univalent vertices of $\Gamma$.  $X$ is the {\em skeleton}
of $\Gamma$.  The {\em degree} of such a
diagram is half the number of trivalent vertices of $\Gamma$.
\end{definition}

\begin{definition} For a compact 1-manifold $X$, the space $\cA^f(X)$ is the
vector space spanned by Jacobi diagrams based on $X$ modulo the
antisymmetry and vertex relations of Equations~\ref{eq:antisymmetry}
and~\ref{eq:vertex}.  $\cA(X)$ is the completion of $\cA^f(X)$ with
respect to the degree.
\end{definition}

\begin{definition}
$\cA^\bc(X)$ is the subspace of $\cA(X)$ spanned by boundary-connected
Jacobi diagrams: diagrams with no connected components that do not
contain a portion of the skeleton $X$.
\end{definition}

By the construction of Chapter~\ref{chap:kontsevich}, $\cA^\bc(X)$ is
isomorphic to the space $\cK_n(X)/\cK_{n+1}(X)$.

To relate these definitions to the definitions in
Section~\ref{sec:diagrams1}, note that $\cA(\uparrow)$ (i.e., $\cA$ of
an oriented interval) is isomorphic to $\cA$.  Recall
that $\cA$ was defined to be the space of Jacobi diagrams with some
ordered, univalent ends, modulo the antisymmetry, Jacobi, and STU
relations.  To pass from a diagram in $\cA$ to a diagram in
$\cA(\uparrow)$, place the univalent vertices on the interval in the
specified order, always attaching from the left side of the interval
(as in Figure~\ref{fig:A-examp}).  The antisymmetry relation in $\cA$
becomes an antisymmetry relation in $\cA(\uparrow)$ and the Jacobi and
STU relations both become vertex relations (of the two different
types) in $\cA(\uparrow)$.  The inverse map from $\cA(\uparrow)$ to
$\cA$ is applying antisymmetry relations so all legs attach from the
left and dropping the interval, remembering the order.  In the future
we will not distinguish between these two spaces.

We will also use versions of the space $\cB$ in this more general
context.  Recall that $\cB$ is the space $\cA$ with the ordering on
legs (and the STU relation) dropped.  We can perform a similar
operation on any interval component of $X$.

\begin{definition}
A {\em Jacobi diagram based on $X\cup Y$}, where $X$ is a 1-manifold
and $Y$ is a set of asterisks $\ast$, is a vertex-oriented
uni-trivalent graph $\Gamma$ with a proper embedding
$\iota:X\rightarrow\Gamma$ up to isotopy with the univalent
vertices that are not in the image of $\partial X$ labelled by
asterisks in $Y$.  In addition to self-loops and multiple edges,
circle components with no vertices are allowed, as long as they are in
the image of $\iota$.
The {\em degree} of $\Gamma$ is half the number of
trivalent plus $Y$-labelled univalent ends.
\end{definition}

\begin{definition}
The space $\cA^f(X\cup Y)$, $X$ and $Y$ as above, is the space of
Jacobi diagrams on $X\cup Y$ modulo antisymmetry relations at each
trivalent vertex and vertex relations around each trivalent vertex.
$\cA(X\cup Y)$ is the completion of $\cA^f(X\cup Y)$ with respect to
degree.
\end{definition}

For $X$ a collection of 1-manifolds and/or asterisks, there is a map
$\chi_x:\cA(\ast_x\cup X) \rightarrow\cA(\uparrow_x\cup X)$, defined
analogously to Definition~\ref{def:chi}: take the average over all
possible ways of ordering the univalent legs labelled by $x$ and
attach them to an oriented interval on the left side.  The proofs that
$\chi$ is a vector space isomorphism also prove that $\chi_x$ is a
vector space isomorphism.

If we want to similarly symmetrize over circle components, we run into
a snag.  There is a natural map from $\cA(\vert\cup X)$ to
$\cA(\bigcircle \cup X)$ given by attaching the two
univalent ends in the diagram that are the images of the endpoints of
the interval $\vert$.  In the simplest case this map is an isomorphism.


\begin{lemma}\lbl{lem:line-circle}
The spaces $\cA(\uparrow)$ and $\cA(\directedcircle)$ are 
isomorphic.
\end{lemma}

\begin{proof} The map from $\cA(\uparrow)$ to $\cA(\directedcircle)$
defined above is clearly surjective.  To see that it is injective, we
need to check that diagrams differing by cyclic permutations are
already equal in $\cA(\uparrow)$.  It suffices to check that two
Jacobi diagrams based on an interval differing by moving a single leg
from the beginning to the end of the interval are equal modulo the
antisymmetry and vertex relations.  Since vacuum diagrams play no part
in this problem, let us assume that there are none.  We can then pick
a good routing of the diagram and consider an associated
$n$-circus.  Now we can expand the first loop on the interval,
pass it over the rest of the circus, and shrink it down again; we
have moved this leg from the beginning to the end of the interval.
See Figure~\ref{fig:line-circle}.

This can also be achieved by a sequence of vertex and antisymmetry
relations without referring to topology by replacing the loops in
Figure~\ref{fig:line-circle} by rounded boxes as in
Definition~\ref{def:rbox} and sweeping the boxes from one end of the
interval to the other, applying vertex and antisymmetry relations
along the way.
\end{proof}

\begin{figure}\label{fig:line-circle}
\[
\figcent{line-circle-1}{0.3}
= \figcent{line-circle-2}{0.3}
= \figcent{line-circle-3}{0.3}
= \figcent{line-circle-4}{0.3}
\]
\caption{A leg in $\cA(\vert)$ can be moved from one end of the
interval to the other.}
\end{figure}

The same proof works to show that $\cA(\uparrow\cup X) \isom
\cA(\directedcircle\cup X)$ where $X$ is a closed 1-manifold; but
the proof does not work (and the statement is not true) if there is
another interval component.  Explicitly,
$\cA(\uparrow\uparrow)\nisom\cA(\uparrow\directedcircle)\isom\cA(\directedcircle\directedcircle)$.
The problem is that as we sweep the loops they get ``caught up'' on the
extra ends of $X$.  So if we want to symmetrize over the legs on circle
components in general, we need to add a new relation.

\begin{definition} In $\cA(\ast_y \cup X)$, {\em link relations} on
$y$ are parametrized by Jacobi diagrams based on $\ast_y \cup X$ in
which one of the $y$-labelled legs is distinguished and marked by
`$*y$'.  The corresponding link relation is the sum of all ways of
attaching the marked leg to all the other legs labelled $y$:
\[
\figcent{link-rel-0}{0.4} \mapsto
\figcent{link-rel-1}{0.4} =
  \figcent{link-rel-2}{0.4}
+ \figcent{link-rel-3}{0.4} + \cdots
+ \figcent{link-rel-4}{0.4}.
\]
\end{definition}

Link relations are the image in $\cA(\ast_y\cup X)$ of moving one
strand from the beginning to end of the interval in
$\cA(\uparrow_y\cup X)$:
\[
   \figcent{A-link-rel-1}{0.4} - \figcent{A-link-rel-2}{0.4}
   = \figcent{A-link-rel-3}{0.4}
   = \figcent{A-link-rel-4}{0.4}.
\]
Symmetrizing the last diagram leaves it in the same form.  So, inside
$\cA(\directedcircle_y\cup X)$, all link relations vanish.
Conversely, if we take the quotient of $\cA(\ast_y\cup X)$ by all link
relations on $y$, we get a space isomorphic to
$\cA(\directedcircle_y\cup X)$.

\begin{definition} $\cA(X \cup Y \cup Y')$, where $X$ is a 1-manifold,
$Y$ is a set of asterisks $\ast$, and $Y'$ is a set of
circled asterisks $\circledast$, is the space of Jacobi diagrams based
on $X \cup Y \cup Y'$ modulo the vertex and antisymmetry relations as
before and, in addition,
link relations on each label in $Y'$.
\end{definition}

\chapter{The Kontsevich integral and its properties}
\lbl{chap:kontsevich}

In this chapter, we will briefly review the {\em Kontsevich integral},
a remarkable invariant of tangles that is a universal
Vassiliev invariant.

\begin{definition}\lbl{def:univ-vassiliev} A {\em universal Vassiliev
invariant} $Z$ is a map $Z: \cK(X) \rightarrow \cA(X)$ so that, for any
$n$-circus $C$ in the complement of an $X$-tangle,
\[
Z(\delta^{(n)}(C)) = \Diag'(C) + \text{higher order terms}
\]
Recall that $\Diag':\cK(X) \rightarrow \cA_n(X)$ is the approximation
to $\Diag:\cK(X) \rightarrow \cK_n(X)/\cK_{n+1}(X)$.
\end{definition}

A universal invariant provides an inverse to the map $m: \cA_n(X)
\rightarrow \cK_n(X)/\cK_{n+1}(X)$ constructed in
Section~\ref{sec:weight-systems}, and so proves that the two spaces
are isomorphic.  Furthermore, such an invariant is universal in the
category-theoretical sense: any finite-type invariant with values in a
vector space factors through $Z$.

We will not provide complete details here; for complete proofs, 
see~\cite{BarNatan:OnVassiliev, LeMurakami:Parallel,Lescop:Kontsevich}.  Our main interest is in
the good properties the Kontsevich integral satisfies under a few
natural topological operations; we will consider
connected sum and connected and disconnected cabling.
Except for the computation in Section~\ref{sec:conn-cable}, which is
due to Thang Le, all material in this chapter is standard.


\section{The Kontsevich integral for braids} We will first define the
Kontsevich integral for braids.  A {\em braid} on $n$ strands for our
purposes is a smooth map from the interval $[0,1]$ to the
configuration space of $n$ points in $\RR^2$.  There is a canonical
framing for braids, since the vector field pointing in (say) the $y$
direction is always normal to the braid.  We can multiply two braids
if the endpoint of one is the starting point of the other, just by
placing one on top of the other and smoothing a little.

The Kontsevich integral for braids is the holonomy of a certain formal
connection on the configuration space of $n$ distinct points in the
plane called the formal Knizhnik-Zamolodchikov (KZ) connection.  It takes
values in the space of chord diagrams appropriate to $n$-strand
braids,
\[
\cA(\underbrace{\vert \cdots \vert}_{\text{$n$ copies}})
   \overset{\text{def}}{=} \cA(n).
\]
Note that $\cA(n)$ has an algebra structure given by vertical composition,
just like braids.  The values actually lie in the subspace of this
space generated by horizontal chords.

Let $c_{ij}$ be a chord connecting the $i$'th and $j$'th strands,
\[
c_{ij} = \figcent{chordij}{0.2}
\]
Let $z_j = x_j + i y_j$ be complex coordinates for the configuration
space of $n$ points.  Then the connection 1-form of the KZ connection is
\[
\mu_n = \frac{1}{2\pi i}\sum_{i < j} c_{ij} d\,\log(z_i - z_j)
      = \frac{1}{2\pi i}\sum_{i < j} c_{ij} \frac{dz_i - dz_j}{z_i-z_j}.
\]

\begin{exercise}
Show that the Knizhnik-Zamolodchikov connection is flat:
\[
d\mu_n + [\mu_n, \mu_n] = 0.
\]
\end{exercise}

The Kontsevich integral of a braid is defined to be the holonomy of this
connection, defined using Chen's iterated integral
construction~\cite{Chen:Iterated}
\[
Z(B) = \sum_{m=0}^{\infty} \int_{0 < t_1 < \cdots < t_m < 1}
         B^*(\mu_n)(t_1) \cdots B^*(\mu_n)(t_m).
\]
This is invariant under isotopy of the braid since the
Knizhnik-Zamolodchikov connection is flat.

More explicitly, we can write this as
\ssp
\begin{equation} Z(K)=\sum_{m=0}^\infty \frac{1}{(2\pi i)^m}
	\mathop{\int}_{t_1<\ldots<t_m}
	\sum_{\substack{\text{applicable pairings} \\
		P=\{(z_i,z'_i)\}}}	D_P
	\bigwedge_{i=1}^{m}\frac{dz_i-dz'_i}{z_i-z'_i}
	\in\cA(n). \label{eq:kontsevich-integral}
\end{equation}
\dsp

In the above equation,
\begin{itemize}
\item an `applicable pairing' is a choice of an unordered
	pair $(z_i,z'_i)$ for every $1\leq i\leq m$, for which
	$(z_i,t_i)$ and $(z'_i,t_i)$ are distinct points on $K$.
\item $D_P$ is the chord diagram naturally associated with $K$ and
      $P$, an appropriate product of the $c_{ij}$'s.
\item Every pairing defines a map $\{t_i\}\mapsto\{(z_i,z'_i)\}$ locally
	around the current values of the $t_i$'s. Use this map to pull
	the $dz_i$'s and $dz'_i$'s to the space $t_1<\dots<t_m$ and
	then integrate the indicated wedge product over that simplex.
\end{itemize}

\section{First properties: Universal and
Grouplike}\label{sec:first-prop} Let us first check that this
invariant is, indeed, a universal finite-type invariant.  The
statement of universality has to be modified slightly for braids,
since an action of an arbitrary circus will take us out of braids and
into general $n$-strand tangles.  We will restrict to {\em horizontal
circuses}: circuses in which each double lasso lies in a horizontal
(constant $t$) plane.  Action on such a circus yields another braid.
We have to show that Kontsevich integral of the resolution of a
horizontal $m$-circus is the corresponding graph plus higher order
terms.

Isotop the braid and the lassos so that all the lassos are very short,
connecting strands close to each other.  Also apply the splitting
relation so that each loop of a lasso encloses a single strand.  When
we resolve a single lasso, we get a difference of two terms differing
only near the crossing.  The KZ connection will be constant on chords
that do not touch the crossing and nearly constant on chords only one
end of which reaches the crossing.  For the alternating difference to
be non-trivial, there must be at least one chord connecting the two
strands involved in the lasso.  Similarly, if we resolve an entire
circus, the integral will be constant or nearly constant if there is
not at least one chord per lasso connecting the two strands of the
lasso.  Thus the lowest degree contribution to $Z(B_D)$ is at least
degree~$m$, and the degree~$m$ term is proportional to the desired
chord diagram.  To check the constant of proportionality, note that we
are integrating $d(log z)/2\pi i$ on a counterclockwise circle around
0; by Cauchy's theorem, the answer is~1.

Next we will show that the Kontsevich integral is {\em grouplike}.
For motivation, consider that the product of two finite type
invariants $v_1$, $v_2$ is again finite type.  The Kontsevich
integral, as the universal finite type invariant, encodes the
information about $v_1$, $v_2$, and $v_1\cdot v_2$; the grouplike
property encodes how the three invariants are related.

\begin{definition}\label{def:box}
For any 1-manifold $X$, the coproduct $\square: \cA(X) \rightarrow \cA(X)
\otimes \cA(X)$ takes a Jacobi diagram $D$ to the sum, over all
partitions $D = D_1 \cup D_2$ of $D$ into two parts in which both
$D_1$ and $D_2$ contain $X$ but are otherwise disjoint, of $D_1
\otimes D_2$.  In other words, assign each connected component of $D
\setminus X$ to the first or second tensor factor and sum over all
possibilities.
\end{definition}

\begin{exercise} Let $v_1,v_2: \cK \rightarrow \CC$ be
Vassiliev invariants of degree $n_1$ with weight systems $w_i:\cA
\rightarrow\CC$.  ($w_i$ is supported on Jacobi diagrams of degree $n_i$.)
Show that $v_1\cdot v_2$ is a Vassiliev invariant of
degree~$n_1 + n_2$ and weight system 
\[
(w_1\otimes w_2) \circ \square: \cA \rightarrow \CC
\]
\end{exercise}

\begin{hint}
The evaulation of $v_1 v_2$ on the resolution of a $k$-circus can be
written as a sum of terms
\[
(v_1\cdot v_2)(\delta^{(k)}(C)) = \sum v_1(\delta^{(l)}(C_1))v_2(\delta^{(k-l)}(C_2))
\]
where $C_1$ and $C_2$
are appropriate circuses of degrees $l$ and $k-l$.
(Think about proving that the product of polynomials
is polynomial.)  For $k = n_1+n_2+1$, at least one of the two factors
vanishes.  For $k=n_1+n_2$, the sum divides the lassos of $C$ into two
subsets of size $n_1$ and $n_2$.  For circuses related to graphs, each
connected component must go into the same subset or the result is 0.
\end{hint}

\begin{proposition}
For any braid $B$, $Z(B)$ is {\em grouplike}:
\[
\square Z(B) = Z(B) \otimes Z(B).
\]
\end{proposition}

\begin{proof}[Proof (Sketch)] This is a general property of the
holonomy of a connection with values in a Hopf algebra in which the
connection form is primitive.  (The usual example is the holonomy of a
Lie-algebra valued connection, which takes values in the the Lie group
or, alternatively, the grouplike elements inside the universal
enveloping algebra.)
\end{proof}

By the general structure theory of Hopf algebras, any time there is a
multiplication on $\cA(X)$ compatible with the comultiplication
$\square$, a grouplike element is an exponential of primitive elements.
Here a {\em primitive} element is an element $x$ such
that $\square(x) = x \otimes 1 + 1 \otimes x$.
In particular, the multiplications in $\cA$, $\cB$, and $\cA(n)$ are
compatible with $\square$.

\section{Cups and caps}\label{sec:cups-caps} In this thesis, we are
not interested in braids, but rather in knots and, more generally,
tangles.  The difference is the presence of critical points, either
minima $\cup$ or maxima $\cap$.  The integral in
Equation~\ref{eq:kontsevich-integral} naturally extends to this
context.  But there is a problem: the integral no longer converges.
For instance, the integral corresponding to a single chord at a maximum
$\figcent{problem-integral}{0.025}$, with the strands are a distance
$d$ apart at the bottom, is
\[
\int_{z=d}^0 \frac{1}{2\pi i}\frac{dz}{z} \cdot
\text{\LARGE$\isolatedchord$}
\]
which diverges logarithmically.  

We can fix the divergence by
terminating the integral when the points are a displacement $\epsilon$
apart and multiplying by the counterterm \(\exp(-\frac{\log\epsilon}{2\pi
i})\isolatedchord\).  It turns out that there is then a finite limit
as $\epsilon\rightarrow 0$. To make the counterterm well-defined, we
need to pick
a branch of $\log\epsilon$ near each critical point of the height
function.  For concreteness, turn the knot near each critical point so
that the tangent at the critical point is along the $x$ axis (so
$\epsilon$ is real) and take the real value of
$\log\,\lvert\epsilon\rvert$ in the counter term.  A knot with the
critical points pinned down like this has a canonical framing, since a
vector field pointing in the $y$ direction is transverse to the knot.
(If we project the knot onto the $x$-$t$ plane, this is the blackboard
framing.)  If we twist a maximum counterclockwise or a minimum
clockwise by $\pi$, the counter term changes by
$\exp(\frac{1}{2}\isolatedchord)$ and the above framing changes by 1
unit, so pinning down the critical points like this is equivalent to a
choice of framing.  See Le and Murakami~\cite{LeMurakami:Universal}
or Lescop~\cite{Lescop:Kontsevich} for details.

For any $X$-tangle $T$, let $\tilde Z(T)$ be the renormalized $\cA(X)$
valued Kontsevich integral as above.  Although $\tilde Z$ converges,
it is not an invariant.  In particular, \(\tilde Z(\wiggle) \ne
\tilde Z(\vert) = \vert.\) Let $\nu^{-1} = \tilde Z(\wiggle)$.

\begin{definition} The Kontsevich integral $Z(T)$ is
$\tilde Z(T)$ multiplied by $\nu$ at each local maximum.
\end{definition}

Now whenever we want to straighten a wiggle $\wiggle$, there is an
extra $\nu$ to cancel the $\nu^{-1}$ from the non-invariance of
$\tilde Z(\wiggle)$.  With these definitions, $Z(\bigcirc) = \nu.$
This element, $\nu$, is what we will compute in
Chapter~\ref{chap:wheels}.

It is straightforward to see that this extension of $Z$ to tangles is
still grouplike and a universal finite-type invariant.

From the discussion above about the framing and the critical points,
we can compute what happens if we change the framing on a component.

\begin{lemma}\label{lem:framing-change}
Given a knot $K$, let $K'$ be the same knot with the framing twisted
by $f$ positive twists.  Then
\[
Z(K') = Z(K) \connect \exp\left(\frac{f}{2}\mlarge\isolatedchord\right).
\]
Similar results hold for changing the framing on a component in a
tangle, where we multiply by $\exp(\frac{f}{2}\isolatedchord)$ on the
affected component.
\qed
\end{lemma}

With only a little more work we can be more precise.  The integrals
contributing to a single chord $\isolatedchord$ come from the term $m=1$ in
Equation~\ref{eq:kontsevich-integral}.  Since all such integrals are of
the form $\frac{1}{2\pi i}\int \frac{dz}{z} = \frac{\log z}{2\pi i}$,
this term can be computed.  Mostly it computes the winding numbers of
strands about themselves.
Careful accounting of what happens at the maxima and minima shows that
the coefficient of $\isolatedchord$ is $1/2$ the number of positive
crossings minus the number of negative crossings.  If we combine this
with the grouplike property of the Kontsevich integral, we get the
following lemma.

\begin{lemma}\label{lem:struts}
The Kontsevich integral of a knot $K$ with framing $f$ relative to the
canonical framing (the framing given by, e.g., a Seifert surface for
$K$) is
\[
Z(K) = \exp_{\connect} \left(\frac{f}{2}\strutn{}{}\right) \connect D
		\in \cB
\]
where $D$ is strutless: has no connected components which are struts
$\strutn{}{}$.\qed
\end{lemma}

\section{Connected sum} \label{sec:connect-sum} The first thing to
notice about our definition of $Z$ is that it is local in the
horizontal plane: each horizontal slice is independent.  Therefore,
for any tangles $T_1$, $T_2$ so that the upper boundary of $T_1$ is
the same as the lower boundary of $T_2$,
\begin{equation}
Z(T_1 \cdot T_2) = Z(T_1) \cdot Z(T_2)
\end{equation}
where the multiplications of tangles on the LHS and web diagrams on
the RHS are both vertical stacking: place the two objects on top of
each other and join the corresponding legs.

\begin{figure}[htbp]
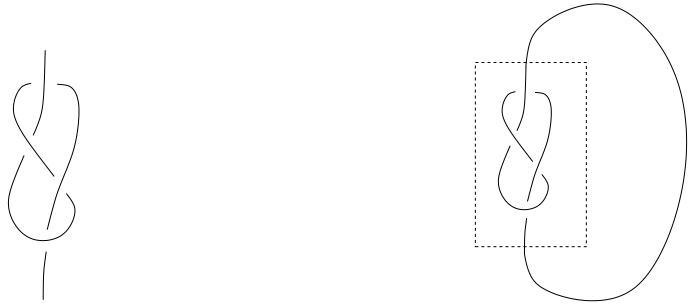
%
\begin{center}%
   {\makebox[.3\linewidth]{\figscaled{figure8-open}{0.4}}}%
\hfil%
   {\makebox[.3\linewidth]{\figscaled{figure8-closed}{0.3}}}%
\end{center}%
\caption{A $(1,1)$ tangle and its closure.}\label{fig:11tangle}%
\end{figure}

We will be particularly interested in the case where the tangles are
$(1,1)$ tangles, with one incoming and one outgoing strand.  These are
closely related to knots: we can turn a $(1,1)$ tangle canonically
into a knot by joining the top to the bottom.  See
Figure~\ref{fig:11tangle} for an example.  What does this do the the
Kontsevich integral $Z$?  By shrinking the tangle, the dotted box in
Figure~\ref{fig:11tangle}, the interactions between the dotted box and
the original strand in the integral will be go to 0 while leaving the
integral within the box unchanged.  The
differences between the invariants of the two are caused by (a)~the
correction introduced at the maximum and (b)~changing the space of
values by closing the strand, which by Lemma~\ref{lem:line-circle} is
isomorphism.  Letting $T$ by the tangle and $K$ be its closure, which
could be a knot or a link with a distinguished component, we find
\begin{equation}\label{eq:11tangle}
Z(K) = \nu \connect Z(T)
\end{equation}
where the two spaces are identified by their isomorphism.

Let $T_1$, $T_2$ be two $(1,1)$ tangles and $K_1$, $K_2$ their
closures.  Then
\begin{align*}
Z(K_1 \connect K_2) &= Z(T_1 \cdot T_2) \connect \nu\\
		    &= Z(T_1) \connect Z(T_2) \connect \nu \\
		    &= Z(K_1) \connect Z(K_2) \connect \nu^{-1}.
\end{align*}

\section{Cables} \label{sec:cables}
\begin{figure}[tbp]
\begin{center}%
\subfigure[Disconnected cable.]%
  {\makebox[.3\linewidth]{\fig{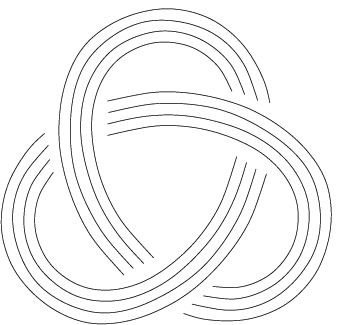}{0.3}%
			\label{fig:disconn-cable-ex}}}%
\subfigure[Connected cable.]%
  {\makebox[.3\linewidth]{\fig{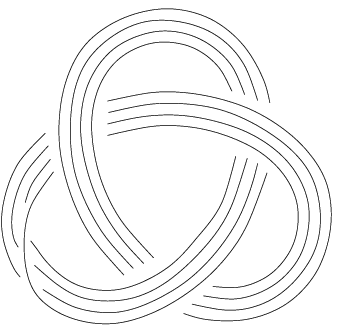}{0.3}%
			\label{fig:conn-cable-ex}}}%
\subfigure[General satellite: the Whitehead double]%
  {\makebox[.3\linewidth]{\fig{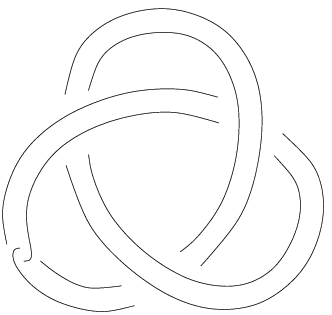}{0.3}%
			\label{fig:satellite-ex}}}%
\end{center}
\caption{Example cablings and satellites of the trefoil knot.}
\label{fig:cable-ex}
\end{figure}

We now consider the operation of {\em cabling}: replacing a knot
with $n$ parallel copies of the same knot, as in
Figure~\ref{fig:cable-ex}.  There are two versions: the
disconnected cabling, as in Figure~\ref{fig:disconn-cable-ex}, ending
up with an $n$-component link, and the connected cabling, as in
Figure~\ref{fig:conn-cable-ex}, in which you add a twist so the
result is again a knot.

Note that this operation is only well-defined for framed knots, and
for connected cabling, you need to specify the extra twist to be
added.  It is natural to include this extra twist in the framing as a
kind of ``rational framing.''

\subsection{General Satellites}
Let us start with some generalities about finite type invariants of
satellites.  A general satellite operation is specified by a link $L$
embedded in a solid torus; the operation is to take a knot $K$, remove
a tubular neighborhood of it, and glue in the solid torus with the
embedded link, obtaining a new link $K_L$.  An example of the result is in
Figure~\ref{fig:satellite-ex}.

What does this operation do to finite-type invariants?  Consider an
$n$-circus in the complement of $K$.  After the satellite
operation, it becomes an $n$-circus in the complement of $K_L$.  By
the Splitting Relation in Section~\ref{sec:weight-systems}, the loops
of the double lassos circling $K$ turn into a sum over the components
of $L$ at that point.

Let us formalize this.

\begin{definition} Let $X, Y$ be compact 1-manifolds and $\phi: X
\rightarrow Y$ a proper map.  Let $D$ be a diagram on $Y$
with a fixed parameterization of $Y$, so that the legs of $D$ have a
well-defined position in $Y$, which are generic with respect to
critical points of $\phi$.  A {\em lift} of $D$ by $\phi$ is a diagram
$D'$ on $X$ with the same internal part on $X$ so that the image of
each leg of $D'$ under $\phi$ is the corresponding leg of $D$.  Each
lift comes with an orientation induced from the local map.
\end{definition}

\begin{definition}\label{def:pullback} Let $X, Y$ be compact
1-manifolds and $\phi: X\rightarrow Y$ a proper map.  The {\em
pullback} $\phi^*: \cA(Y)\rightarrow \cA(X)$ on a diagram $D \in
\cA(Y)$ is the sum over all lifts of $D$ from $Y$ to $X$.
\end{definition}

$\phi^*$ is well-defined (i.e., does not depend on the position of the
legs of $D$ and descends modulo the relation) and is invariant under
homotopy of $\phi$.

\begin{proposition}[Kuperberg~\cite{Kuperberg:Invertibility}] Let $X,
Y$ be 1-manifolds, let $\eta: Y\rightarrow N(X)$ be a pattern for
satelliting, and let $\iota: X\rightarrow\RR^3$ be a link/tangle with
an $n$-circus $C$ in its complement.  Then for any universal invariant
$Z$, $Z(\delta_C(X_\eta)) = \phi^*(\Diag(C)) + \text{h.o.t}$, where
$\phi$ is the composition of $\eta$ and the retraction
$N(X)\rightarrow X$.
\end{proposition}

\begin{proof}
By the definition of a universal finite type invariant $Z$,
\[
Z(\delta_C(X_\eta)) = \Diag(C_\eta) + \text{h.o.t.},
\]
where $C_\eta$ is the $n$-circus $C$ in the complement of $X_\eta$.  By
the Splitting relation, this is $\phi^*(\Diag(C))$.
\end{proof}

\begin{corollary}\label{cor:satellite}
The map $\eta^*: \cA(X) \rightarrow \cA(Y)$ is $\phi^*$ plus higher
order terms (terms that increase the degree of a diagram).
\end{corollary}

The nice property of the two cabling operations is that the above
formula becomes exact: the ``higher order terms'' vanish.

\begin{definition} The operations
\[
\Delta^x_{x_1\dots x_n}: \cA(\vert_x\cup X)\rightarrow
	\cA(\vert_{x_1}\cup\cdots\cup\vert_{x_n}\cup X)
\]
or
\[
\Delta^x_{x_1\dots x_n}: \cA(\bigcircle_x\cup X)\rightarrow
	\cA(\bigcircle_{x_1}\cup\cdots\cup\bigcircle_{x_n}\cup\,X)
\]
are the pullback of the $n$-fold disconnected cover of the component
labelled $x$.  When we do not care about the labels on the result, an
alternate notation is $\Delta^{(n)}_x$.
\end{definition}

\begin{proposition}[Le and Murakami~\cite{LeMurakami:Parallel}]
\label{prop:disconn-cable}
Let $L$ be a framed tangle and $\DCable_x^n(L)$ be the $n$-fold
disconnected cable of $L$ along a knot component $x$ or interval
component with one upper and one lower boundary.  Then
\[
Z(\DCable_x^n(L)) = \Delta^{(n)}_x(Z(L)).
\]
\end{proposition}

Proof in Section~\ref{sec:disconn-cable}.

\begin{definition} The operation
\[
\psi^{(n)}_x: \cA(\bigcirc_x\cup X)\rightarrow\cA(\bigcirc_x\cup X)
\]
is the pullback of the $n$-fold connected
cover of the circle labelled $x$.
\end{definition}

\begin{proposition}[T.~Le]
\label{prop:conn-cable}
Let $L$ be a framed tangle and $\CCable_x^n(L)$ be the $n$-fold
connected cable of $L$ along a knot component $x$ as in
Figure~\ref{fig:conn-cable-ex}.  Then
\[
Z(\CCable_x^n(L)) =
	\psi^{(n)}_x(Z(L) \connect_x \exp(\frac{1}{2n}\isolatedchord))
\]
for an appropriate choice of framing on $\CCable_x^n(L)$.
\end{proposition}

Proof in Section~\ref{sec:conn-cable}.

One reason to introduce symmetrized diagrams is that the
operations $\Delta$ and $\psi$ above become very simple in $\cB$.

\begin{lemma}\label{lem:delta} The map 
\[
\Delta^x_{x_1\dots x_n}:\cA(\ast_x \cup X) \rightarrow
	\cA(\ast_{x_1} \cup\dots\cup\ast_{x_n} \cup X)
\]
is the sum over all ways of replacing each $x$ leg by one of the
$x_i$.\qed
\end{lemma}

\begin{remark}
$\Delta$ is similar to a coassociative, cocommutative coproduct in a
coalgebra, except that it does not take values in $\cA \otimes \cA$.
Do not confuse it with the coproduct $\square$ of
Definition~\ref{def:box}, which is an honest coproduct.
\end{remark}

The operation $\Delta$ in Lemma~\ref{lem:delta} is analogous to a
change of variables $x \mapsto x_1 + \dots + x_n$ for ordinary
functions $f(x)$.  We will use a suggestive notation: a leg labelled
by a linear combination of variables means the sum over all ways of
picking a variable from the linear combination.  If $D(x)$ is a
diagram with some legs labelled $x$, $\Delta^{(n)}(D(x)) = D(x_1 +
\dots + x_n)$ is the diagram with the same legs labelled $x_1 + \dots
+ x_n$.

\begin{lemma}\label{lem:psi} The map
\[
\psi^{(n)}_x: \cA(\circledast_x \cup X) \rightarrow
	\cA(\circledast_x \cup X)
\]
is multiplication by $n^k$ on diagrams with $k$ legs labelled $x$.
\qed
\end{lemma}

This operation is related to the change of variables $x \mapsto nx$.

\subsection{Disconnected cabling}\label{sec:disconn-cable}
Here we sketch the proof of Proposition~\ref{prop:disconn-cable}.  For
full details, please see Le and Murakami~\cite{LeMurakami:Parallel}.

First consider the disconnected cabling of a
braid.  Take one strand of a braid and replace it by a
number of copies parallel with a small displacement $\epsilon$.
Because the strands are parallel, there are no new chords between them
($z_1 - z_2$ is constant); since the new strands are close to the
original strand, thep
integrals for chords between the doubled strands and another strand
are close to chords for the original strand, with equality in the
limit as $\epsilon \rightarrow 0$.  See
Lescop~\cite{Lescop:Kontsevich} for the analytic details.  Hence for
parallel cabling of a braid Proposition~\ref{prop:disconn-cable}
holds.

When we switch to tangles, there may be cups and caps.  The Kontsevich
integral of a doubled maximum $Z(\doublecap)$ will not be the double
of the 
Kontsevich integral of a single maximum $Z(\cap)$.  

\begin{figure}
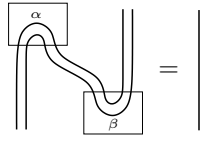
\label{fig:wiggle2}
\[
\figcent{wiggle2}{0.2} = \figcent{id2}{0.2}
\]
\caption{A doubled wiggle.}
\end{figure}

Let $\alpha = \lim_{\epsilon\rightarrow 0} Z(\doublecap_\epsilon)$
where $\doublecap_\epsilon$ is a doubled maximum with the distances
between the strands at the bottom $\epsilon$, $1$, and $\epsilon$ in
order.  (There are naively two logarithmic divergences in this
integral which fortuitously cancel.) Note that $\alpha$
includes renormalization of $\nu$ on each strand added so that
$Z(\wiggle) = Z(\vert)$.  Similarly we get another element
$\lim_{\epsilon\rightarrow 0}Z(\doublecup_\epsilon) =
\beta\in\cA(\uparrow\,\uparrow)$ at minima.  Although we can say
little about $\alpha$ or $\beta$ by themselves,%
\footnote{Le and Murakami~\cite{LeMurakami:Parallel} have analyzed
contexts in which you can compute $\alpha$ and $\beta$ and, in
particular, when $\alpha$ and $\beta$ are the doubles of a maximum and
minimum.}
we can compute their
product.  
Consider a doubled wiggle, as in Figure~\ref{fig:wiggle2}.
It can be decomposed into three pieces: a doubled maximum, a doubled
minimum, and the doubling of the braid whose closure is $\wiggle$.
The doubled ends of the Kontsevich integral of the doubling of the
braid can slide over any diagram which, like $\alpha$ or $\beta$,
lives on two strands, so we can collect $\alpha$ and $\beta$ at the
beginning.  We then recognize the remaining integral as the double of
the naive integral which gave us $\nu^{-1}$ above.  Since this whole
$(2,2)$ tangle can be straightened to give the identity, we find
\begin{align*}
\alpha \cdot \beta \cdot \Delta\nu^{-1} &= 1 \\
\intertext{or}
\alpha\cdot\beta &= \Delta\nu
\end{align*}
where the multiplication is the natural multiplication on
$\cA(\uparrow\uparrow)$.

Now consider doubling a knot component or a component with one upper
and one lower end.  There will be the same number
$k$ of maxima and minima, so the net error from a pure doubling will
be some product of $k$ $\alpha$'s and $k$ $\beta$'s.  As before, we
can slide the $\beta$'s through the rest of the integral
until they are next to the $\alpha$'s; then each pair cancel to give
$\Delta\nu$, which is exactly the double of the original
renormalization, as desired.
(There is one difficulty we are brushing under the rug: there are
actually two versions of $\alpha$ and $\beta$, depending on which
direction you traverse the maximum, and several different products can
occur.  To show that all pairs cancel, some more complicated version
of Figure~\ref{fig:wiggle2} need to be considered.)

So Proposition~\ref{prop:disconn-cable} is true for doubling a knot
component. 
By iterating this argument, the same thing is true for an $n$-fold
disconnected satellite.\qed

\subsection{Connected cabling}\label{sec:conn-cable} Now we consider
the case of a connected cabling.  The difference between the connected
cabling and the
disconnected cabling above is the extra $1/n$ twist $T_n$ inserted at one
point:
\[
T_n = \figcent{1n-twist}{0.25}.
\]
By isotopy can assume that this twist occurs in a horizontal slice where
all the other strands are vertical.  We can
apply Proposition~\ref{prop:disconn-cable} on the $(n,n)$ ``tangle''
obtained by excising $T_n$.  (This object is not properly a tangle,
since there is a little piece cut out of it.  But we can still compute
its Kontsevich integral by Proposition~\ref{prop:disconn-cable}.)
To complete the computation, we need to
compute $a := Z(T_n)$.

Repeating $T_n$ $n$ times, we get a full twist which we can
compute using Lemma~\ref{lem:framing-change} and the disconnected cabling
of the previous section:
\[
Z\left(\figcent{n-twists}{0.25}\right)
	= Z\left(\figcent{fulltwists}{0.25}\right)
	= \Delta^{(n)}\left(\exp(\frac{1}{2}\isolatedchord)\right)
		\exp(-\frac{1}{2}\isolatedchord)^{\otimes n} =: b.
\]
The notation $\exp(-\frac{1}{2}\isolatedchord)^{\otimes n}$ means $n$
copies of the framing change element
$\exp(-\frac{1}{2}\isolatedchord)$, one on on 
each of the $n$ strands.

The $n$ copies of $T_n$ that appear are not quite the same: they
differ by cyclic permutations of the strands.
If we arrange the $n$ strands at the top and bottom of $T_n$
to be at the vertices of a regular $n$-gon, the strands are
symmetric and $a^n = b$ or
\[
a = b^{\frac1n} = \Delta^{(n)}\left(\exp(\frac{1}{2n}\isolatedchord)\right)
		\exp(-\frac{1}{2n}\isolatedchord)^{\otimes n}.
\]
More generally, we will need to conjugate $T_n$ by some element $C$ to get
the strands symmetric; this implies that
\[
a = c b^{\frac{1}{n}} S(c^{-1})
\]
where $c = Z(C)$ and $S$ is the automorphism of $\cA(\uparrow_{x_1}
\dots\uparrow_{x_n})$ which rotates the strands by $x_i \mapsto x_{i-1}$.

\begin{lemma}\label{lem:n-twist}
$Z(T_n) = c \Delta^{(n)}\left(\exp(\frac{1}{2n}\isolatedchord)\right)
		\exp(-\frac{1}{2n}\isolatedchord)^{\otimes n}
S(c^{-1})$ for some $c \in \cA(\uparrow_{x_1} \dots \uparrow_{x_n})$.
\qed
\end{lemma}

This lemma can also be proved without using the specific geometry of
the Kontsevich integral in an argument due to D.~Bar-Natan, but we
will not give the argument here.



\begin{proof}[Proof of Proposition~\ref{prop:conn-cable}] By the above
computations, the invariant of the connected cable of a knot $K$ is
$\Delta^{(n)}(Z(K))$, multiplied by $Z(T_n)$, and closed up with a
twist.  The conjugating elements $c$ and $S(c^{-1})$ can be swept
through the knot and cancel each other.  The factor
$\Delta^{(n)}(\exp(\isolatedchord/2n))$ in $a$ can be combined with $Z(K)$
so that we apply $\Delta^{(n)}$ to $Z(K) \connect
\exp(\isolatedchord/2n)$.  The twisted closure turns $\Delta^{(n)}$ into
$\psi^{(n)}$.  The remaining $n$ factors of $\exp(-\isolatedchord/2n)$ in
$a$ can be slid around the knot and combined to give
\[
Z(\CCable^n(K)) = \psi^{(n)}\left(Z(K)\connect\exp(\frac{1}{2n}\isolatedchord)\right)
			\connect \exp(-\frac{1}{2}\isolatedchord).
\]
The last factor is absorbed in the ``appropriate choice of
framing'' in the statement of the proposition.
\end{proof}.

\chapter{Wheeling}\label{chap:wheeling}

Now that we have reviewed the basic theory of Vassiliev invariants,
the proof of Theorem~\ref{thm:wheeling} is relatively straightforward.
In Section~\ref{sec:phi} we will show how to interpret the Hopf link
$\OpenHopf$ of Figure~\ref{fig:One} as a map $\Phi: \cB\rightarrow\cA$
by taking its
Kontsevich integral in $\cA(\vert \circledast)$ and gluing the legs on
the bead to the diagram in $\cB$.
In Section~\ref{sec:multiplicativity} we will see how the equation
``$1+1=2$'' of Figure~\ref{fig:OneOneTwo} implies that $\Phi$ is
multiplicative.  Briefly, we interpret both sides of the
equation as maps from $\cB \otimes \cB$ to $\cA$.  The ``$1+1$'' side
takes $(X, Y)$ to $\Phi(X) \connect \Phi(Y)$ and the ``$2$'' side
takes $(X, Y)$ to $\Phi(X \udot Y)$.  The fact that the two links are
isotopic implies that the two maps are equal.  In
Section~\ref{sec:phi0} we relate this map $\Phi$ to the
wheeling map $\Upsilon$ of Section~\ref{sec:upsilon} by showing that
the lowest degree term
$\Phi_0$ of $\Phi$
is also multiplicative and is the same as the Duflo map with possibly
different coefficients of the wheels $\omega_n$.  To fix the coefficients, it
suffices to do a small computation for some Lie algebra; we do the
computation for $\sltwo$ in Appendix~\ref{app:coefficients}.

\section{The map $\Phi$}\label{sec:phi}
We start by defining a kind of inner product on the space $\cB$.

\begin{definition}\label{def:diagram-inner} For a diagrams $C, C' \in
\cB$ so that $C$ has no struts, the {\em inner product} of $C$ and $C'$ is
\ssp
\[ \cA(\emptyset) \contains \langle C, C' \rangle =\begin{cases}
    \parbox{2.2in}{
      the sum of all ways of gluing all the legs of $C$ to all
      legs of $D$
    }\quad & \parbox{1.7in}{if $C$ and $C'$ have
			    the same number of legs} \\[3pt]
    0 & \text{otherwise}
  \end{cases}
\]
\dsp
We will sometimes want to fix $C$ and consider $\langle
C\comma\cdot\rangle$ as a map from
$\cB$ to $\cA(\emptyset)$; we will denote this map $\iota(C)$.  This
definition works equally well in the presence of other skeleton
components or to glue several components.  We
will use subscripts to indicate which ends are glued.
\end{definition}

As in Definition~\ref{def:diagram-differential}, the restriction that
$C$ not have struts is to guarantee convergence and avoid closed
circles.

There are two dualities relating
$\langle\cdot\comma\cdot\rangle$ with other operations we have
defined.

\begin{lemma}\label{lem:mult-comult-dual}
Multiplication and comultiplication in $\cB$ are dual
in the sense that
\[
\langle C, D_1 \udot D_2 \rangle =
	\langle \Delta_{xy} C, (D_1)_x \otimes (D_2)_y\rangle_{xy}.
\]
Similar statements hold in the presence of other ends.
\end{lemma}
\begin{proof} The glued diagrams are the same on the two sides; we
either combine the legs of $D_1$ and $D_2$ into one set and then glue
with $C$, or we split the legs of $C$ into two pieces which are then
glued with $D_1$ and $D_2$. (Note that there are no combinatorial
factors to worry about: in both cases, we take the sum over all
possibilities.)  
\end{proof}

\begin{lemma}\label{lem:mult-diffop-dual} Multiplication by a diagram
$B \in\cB$ and applying $B$ as a diagrammatic differential operator
are adjoint in the sense that
\[
\langle A \udot B, C\rangle = \langle A, \partial_B(C) \rangle.
\]
\end{lemma}
\begin{proof}
As before, the diagrams are the same on both sides.
\end{proof}

The map $\Phi$ of the paper is constructed from the bead on a wire in
Figure~\ref{fig:One}.  We start with its Kontsevich integral:
\[
Z(\openhopf{x}{z}) \in \cA(\uparrow_z, \directedcircle_x).
\]
We then symmetrize on legs attached to the bead $x$ as explained in
Section~\ref{sec:diagrams2}:
\[
\chi^{-1}_x Z(\openhopf{x}{z}) \in \cA(\uparrow_z, \circledast_x).
\]
Finally, we use the inner product operation along the legs
$x$ to get a map from $\cB$ to $\cA$:
\[
\Phi = \iota_x \chi^{-1}_x Z(\openhopf{x}{z}): \cB \rightarrow \cA
\]
In this last step, there are two things we have to check.  First, we
must see that $\chi^{-1}Z(\OpenHopf)$ has no struts.
By Lemma~\ref{lem:struts}, this follows from the fact that we took the
bead with the zero framing.  Second,
we need to check that the inner product descends modulo the link
relations on $x$ in $\cA(\uparrow_z,\circledast_x)$.

\begin{lemma}\label{lem:glue-welldef}
The inner product $\langle\cdot\comma\cdot\rangle_x:\cA(\ast_x\cup X)
\otimes \cA(\ast_x)\rightarrow\cA(X)$ descends to a map
$\langle\cdot\comma\cdot\rangle_x:\cA(\circledast_x\cup X)
\otimes \cA(\ast_x) \rightarrow \cA(X)$.
\end{lemma}

\begin{figure}[htbp]
\[
\mathcenter{\fig{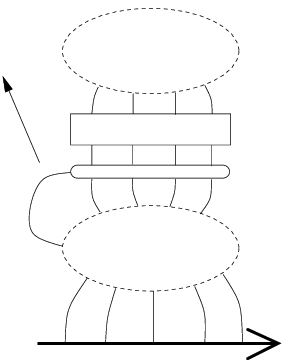}{0.28}}=
\mathcenter{\fig{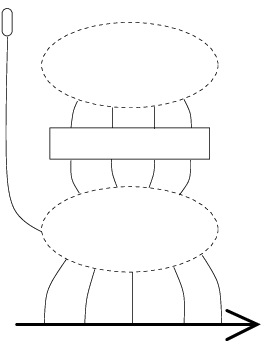}{0.28}}= 0.
\]
\caption{The proof that $\Phi(D)$ is well-defined modulo link
relations on $Z(\OpenHopf)$: link relations in $Z(\OpenHopf)$ can be
slid over $D$.}%
\label{fig:phi-welldef}
\end{figure}

\begin{proof} We use a sliding argument similar to the one in the
proof of Lemma~\ref{lem:line-circle}.  Link relations in
$\cA(\circledast_x\udot X)$ can be slid over diagrams in
$\cA(\ast_x)$, as shown in Figure~\ref{fig:phi-welldef}.
\end{proof}

\section{Multiplicativity of $\Phi$}\label{sec:multiplicativity} We
now come to the key lemma in the proof of Theorem~\ref{thm:wheeling}.

\begin{lemma}\label{lem:phi-multiplicative} The map $\Phi: \cB
\rightarrow \cA$ is an algebra map.
\end{lemma}

\begin{proof}
As advertised, we use the equality of links ``$1+1=2$''.  Let
us see what this equality of links says about the Kontsevich integral
of the Hopf link.  On the ``$1+1$'' side, we see the connected sum of
two open Hopf links; by Section~\ref{sec:connect-sum}, the invariant
of the connected sum is the connected sum of the invariants.  To write
this conveniently, let $H(z;x)$ be $Z(\OpenHopf) \in\cA(\uparrow_z,
\circledast_x)$, with the wire labelled by $z$ and the bead labelled
by $x$.  Then
\[
Z(\OpenHopf \connect \OpenHopf) = H(z; x_1) \connect_z H(z; x_2)
  \in \cA(\uparrow_z, \circledast_{x_1}, \circledast_{x_2}).
\]
On the ``$2$'' side, we see the disconnected cable of a Hopf link.  By
Section~\ref{sec:disconn-cable}, this becomes the coproduct $\Delta$:
\[
Z(\DCable(\OpenHopf)) = \Delta^x_{x_1x_2}H(z;x).
  \in \cA(\uparrow_z, \circledast_{x_1}, \circledast_{x_2}).
\]
Since the two tangles are isotopic, we have
\[
H(z; x_1, x_2) \overset{\text{def}}{=}
	H(z;x_1) \connect_z H(z;x_2) = \Delta^x_{x_1x_2}H(z;x)
  \in \cA(\uparrow_z, \circledast_{x_1}, \circledast_{x_2}).
\]

Now consider the map
\[
\Xi = \iota_{x_1}\iota_{x_2}H(z; x_1, x_2): \cB \otimes \cB \rightarrow \cA;
\]
in other words, in $\Xi(D_1 \otimes D_2)$ glue the $x_1$ and $x_2$
legs of $H(z;x_1,x_2)$ to $D_1$ and $D_2$ respectively.  This descends
modulo the two different link relations in $\cA(\uparrow, \circledast,
\circledast)$ by the argument of Figure~\ref{fig:phi-welldef}, applied
to $D_1$ and $D_2$ separately.  We have two different expressions for
this map from the two different expressions for $H(z; x_1, x_2)$.  On
the ``$1+1$'' side, the gluing does not interact with the connected
sum and we have
\[
\Xi(D_1, D_2) = \Phi(D_1) \connect \Phi(D_2).
\]
See Figure~\ref{fig:glue-1+1}.
For the ``$2$'' side, we use Lemma~\ref{lem:mult-comult-dual} to see that
\begin{align}
\Xi(D_1, D_2) &= \langle \Delta^x H(z; x), D_1 \otimes D_2 \rangle \\
	      &= \langle H(z; x), D_1 \udot D_2 \rangle \\
	      &= \Phi(D_1 \udot D_2).
\end{align}
See Figure~\ref{fig:glue-2}.

\ssp
\begin{figure}[tbp]
\begin{center}
\fig{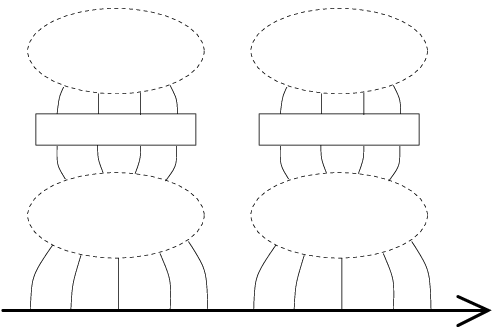}{0.28}
\end{center}
\caption{Gluing $Z(\OpenHopf\connect\OpenHopf)$ to $D_1 \otimes
           D_2$.}\label{fig:glue-1+1}%
\end{figure}
\begin{figure}[tbp]
\begin{center}
\[
\figcent{glue-2a}{0.28} = \figcent{glue-2b}{0.28}
\]
\end{center}
\caption{Gluing $Z(\DCable\OpenHopf)$ to $D_1 \otimes D_2$ in two
equivalent ways.}%
	\label{fig:glue-2}%
\end{figure}
\dsp

Combining the two, we find
\[
\Phi(D_1) \connect \Phi(D_2) = \Phi(D_1 \udot D_2).
\]
\end{proof}

\section{Mapping degrees and the Duflo isomorphism}\label{sec:phi0}
We have successfully constructed a multiplicative map from $\cB$ to
$\cA$.  We will see later (see Section~\ref{sec:Hopf}) that this map $\Phi$
is the Duflo map, but we cannot yet see this.  Instead we will
consider the lowest degree term $\Phi_0$ of $\Phi$.

\begin{definition}\label{def:mapping-degree}
The {\em mapping degree} of a diagram $D \in \cA(\uparrow_z,
\circledast_x)$ with respect to $x$ is the amount $\iota_x D: \cB
\rightarrow \cA$ shifts the degree.  Explicitly, it is the degree
of~$D$ minus the number of $x$ legs of~$D$.
\end{definition}

Since there are no $x$-$x$ struts in $H(z;x)$, every $x$ leg of $H$
must be attached to another vertex (either internal or on the interval
$z$).  Furthermore, if two $x$ legs are attached to the same interal
vertex, the diagram vanishes by antisymmetry.  Therefore there are at
least as many other vertices as $x$ legs in $H$ and the mapping degree
is $\ge 0$.

\begin{definition}\label{def:phi0}
$H_0(z;x)$ is the part of $H(z;x)$ of mapping degree 0 with respect to
$x$.  $\Phi_0$ is $\iota_x H_0(z;x)$.
\end{definition}

$\Phi_0$ is still multiplicative, since the multiplications in $\cA$
and $\cB$ both preserve degrees.  (For homogeneous diagrams $D_1$ and
$D_2$ of degrees $n_1$ and $n_2$, $\Phi_0(D_1 \udot D_2)$ is the piece
of $\Phi(D_1 \udot D_2)$ of degree $n_1+n_2$ and likewise for
$\Phi_0(D_1) \connect \Phi_0(D_2)$.)

\begin{figure}\label{fig:wheels-struts}
\begin{center}
\fig{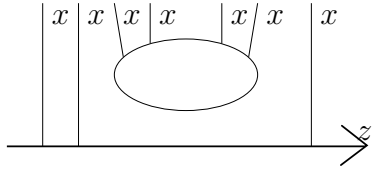}{0.75}
\end{center}
\caption{The only diagrams in $\cA(\uparrow_z, \circledast_x)$ of
mapping degree 0 with respect to $x$ are wheels and struts.}
\end{figure}

The diagrams that appear in $H_0$ are very restricted, since every
vertex that is not an $x$ leg must connect to an $x$ leg; since these
vertices are trivalent, the other two incident edges form a
1-manifold.  The possible diagrams are $x$ wheels and $x-z$ struts, as
shown in Figure~\ref{fig:wheels-struts}.
The linking number between the bead and the wire in the link
$\OpenHopf$ is 1, so the coefficient of the $x$-$z$ strut is 1.
Combined with the fact that the Kontsevich integral is grouplike, we
find that 
\begin{align*}
H_0(z;x) &= \exp(\strutv{x}{z}) \udot \exp(\Omega') \\
\Omega' &= \sum_n a_{2n}\omega_{2n}
\end{align*}
for some coefficients $a_n$.  Note that the right hand side is
written in $\cA(\uparrow, \ast)$ (with a strange mixed product), since
there is no algebra structure on $\cA(\uparrow, \circledast)$.

By the following lemma, we now have a multiplicative map very similar
to our desired map~$\Upsilon$.

\begin{lemma}
$\Phi_0 = \chi \circ \partial_\Omega'$.
\end{lemma}

\begin{proof}
Using Lemma~\ref{lem:mult-diffop-dual} and noting that
gluing with $\exp(\strutv{x}{z})$ takes the legs of a diagram in $\cB$
and averages over all ways of ordering them, as in the definition of
$\chi$, we see that
\[
\Phi_0(D) = \langle \exp(\strutv{x}{z} \udot \Omega', D \rangle
	  = \langle \exp(\strutv{x}{z}, \partial_\Omega'(D) \rangle
	  = \chi(\partial_\Omega'(D)).
\]
\end{proof}

To see that $\Phi_0 = \Upsilon$, we only need to check that $\Omega =
\Omega'$.  This is proved in Proposition~\ref{prop:coefficients}.
\chapter{Wheels and the Kontsevich integral of the unknot}
\label{chap:wheels}
\section{Results}\label{sec:wheels-results}
In this chapter, we will give some explicit computations of the
Kontsevich integral to all orders for some specific knots and links.
We begin by reviewing our principal results.  The first one was
already stated in the introduction.

\setcounter{temp}{\value{theorem}}
\setcounter{theorem}{2}
\begin{theorem}[Wheels; joint with T.~Le]
The Kontsevich integral of the unknot is
\[
Z(\bigcircle) = \nu = \Omega \in \cB.
\]
\end{theorem}
\setcounter{theorem}{\value{temp}}

Recall that $\Omega$ is the ``wheels'' element from the introduction.

By changing the framing on the unknot and cabling it, we can construct
a Hopf link.  Using the results of Section~\ref{sec:conn-cable} and
the value of $Z(\bigcirc)$, we can compute the invariant of the Hopf
link from the invariant of the unknot.
There are several good formulas for the answer.

\begin{theorem}\label{thm:hopf-link}
The framed Kontsevich integral of the Hopf link can be expressed in
the following equivalent ways:
\ssp\begin{align*}
Z(\hopf{x}{y}) &=
  \begin{cases}
    \Upsilon_x\circ\Upsilon_y(\exp(\strutn{y}{x})))\cdot \Vac \\
    \Upsilon_x(\exp_{\udot}(\strutn{y}{x})\Omega_x) \cdot \Vac
  \end{cases} \\
Z(\openhopf{x}{y}) &= \exp(\strutn{y}{x})\udot\Omega_y,
\dsp\end{align*}
for some elements $\Vac \in \cA(\emptyset)$.
\end{theorem}

In the last expression, $\openhopf{x}{y}$ is the $(1,1)$ tangle whose
closure is the Hopf link, with the bead labelled by $y$ and the wire
labelled by $x$.
From this last equality in Theorem~\ref{thm:hopf-link}, we
can see exactly the map $\Phi$ from Chapter~\ref{chap:wheeling}.

\begin{corollary}
$\Phi = \Phi_0 = \chi \circ \partial_\Omega.$
\end{corollary}

\section{Useful facts}

The element $\Omega \in \cB$ is a very remarkable element.  We will
need one of its nice properties for the proof of
the Wheels theorem.  Start from the basic equality proved in the
Wheeling theorem,
\[
\Delta^x_{x_1x_2}H_0(z;x) = H_0(z;x_1) \connect_z H_0(z;x_2)
  \in \cA(\uparrow_z, \circledast_{x_1}, \circledast_{x_2})
\quad \text{where} \quad H_0(z;x) = \Omega_x \exp(\strutn{x}{z}).
\]
Now consider dropping the strand $z$, i.e., mapping all diagrams with
a $z$ vertex to 0.  (Knot-theoretically, this corresponds to dropping
the central strand in the equation ``$1+1=2$''.)  We find
\begin{equation}\label{eq:delta-omega}
\Delta \Omega = \Omega \tensor \Omega \in \cA(\circledast \circledast).
\end{equation}
Note that this equality is not true inside $\cA(\ast\ast)$.

\begin{lemma}[Pseudo-linearity of $\log\Omega$]\label{lem:pseudo-linear}
For any $D \in \cB$,
\[
\partial_D(\Omega) = \partial_D\Omega\rvert_0 \udot \Omega
	= \langle D, \Omega \rangle \Omega.
\]
\end{lemma}
\begin{proof}
\[ \partial_D(\Omega)_x = \langle D_y, \Omega_{x+y}\rangle_y
		        = \langle D_y, \Omega_x \Omega_y\rangle_y
		        = \langle D_y, \Omega_y\rangle_y \Omega_x.
\]
In the second equality, we use Equation~\ref{eq:delta-omega}.  This is
allowed, since the contraction descends to $\cA(\circledast
\circledast) \isom \cA(\circledast \uparrow)$ by
the argument of Lemma~\ref{lem:glue-welldef}.  Note that it is crucial
that $D$ is invariant for this argument.
\end{proof}

\begin{remark} Compare this lemma with standard calculus: if $D$ is
any differential operator and $f$ is a linear function, then
$De^f=(Df)(0)e^f$.  The prefix ``pseudo'' is written above because
Lemma~\ref{lem:pseudo-linear} does not hold for every $D$, but
only for $x$-invariant $D$'s.  (E.g., if $D$ were in $\cA(\ast\ast)$
rather than $\cA(\ast)$, the lemma would not be true.)
\end{remark}

Although we are interested in knots and links in $S^3$ in this thesis,
for which the appropriate space of diagrams $\cA^\bc$ from
Section~\ref{sec:weight-systems} is boundary connected, vacuum
diagrams (elements of $\cA(\emptyset)$)
appear at various points.  Notably, the wheeling map $\Upsilon$ does
not preserve the subspace of boundary connected diagrams.
Although the resulting vacuum
components can be computed explicitly,%
\footnote{D.~Bar-Natan and
R.~Lawrence~\cite{Bar-NatanLawrence:RationalSurgery} have done these
computations}
they are
almost always irrelevant for us and it would just
complicate the formulas to keep track of them.  To avoid this, we will
introduce the {\em boundary-connected projection} $\pi^\bc: \cA
\rightarrow \cA^\bc$ which maps any diagram containing vacuum
components to 0 and is otherwise the identity.  Note
that $\pi^\bc$ is multiplicative.  There are similar projections, which we
also be called $\pi^\bc$, for other spaces $\cA(X)$.

If we compose Lemma~\ref{lem:pseudo-linear} with $\pi^\bc$,
we find
\begin{equation}\label{eq:pseudo-linear-2}
\pi^\bc \partial_D \Omega = \begin{cases}
    \Omega	& \text{$D$ is the empty diagram} \\
    0	& \text{otherwise}
  \end{cases}.
\end{equation}

\section{Coiling the unknot}\label{sec:coiling}

The basic equation we will use to identify $\nu$ is
``$n\cdot 0 = 0$'' from the introduction: the $n$-fold connected cable
of the unknot is the unknot with a new framing.  From
Section~\ref{sec:conn-cable}, this implies that
\begin{equation}\label{eq:n00}
\psi^{(n)}(\nu \connect \exp_{\connect}(\frac{1}{2n}\strutn{}{}))
	 = \nu \connect \exp_{\connect}(\frac{n}{2}\strutn{}{}).
\end{equation}
This equation is true for all $n \in \ZZ$, $n > 0$.  In each degree,
each side is a Laurent polynomial in $n$ of bounded degree; therefore,
the two sides are equal as Laurent polynomials.  The RHS is a
polynomial in $n$, so both sides are polynomials (i.e., have no
negative powers of $n$.)  Let us evaluate both sides at $n=0$.  On the
RHS, we get just $\nu$.  For the LHS, recall how $\psi^{(n)}$ acts in
the space $\cB$: it multiplies a diagram with $k$ legs by $n^k$.  Each
strut appearing in the product on the LHS contributes a factor of
$1/n$; each leg in the symmetrized result gives a factor of $n$.  To
bound how many legs can appear in a symmetrized product, we use the following
lemma.

\begin{lemma}\label{lem:wheels-lemma} For any elements $x_1, \dots,
x_k \in \cA(\uparrow)$ with at least one leg on the interval $\uparrow$,
$\chi^{-1}(x_1 \connect \cdots \connect x_k) \in \cB$ has at least $k$
legs.
\end{lemma}

\begin{proof} First note that any vacuum diagrams that appear in the
$x_i$'s pass through unchanged to the result; let us assume that there
are none, so that we can use the vacuum projection $\pi^\bc$ without
changing the result.  By Theorem~\ref{thm:wheeling},
\[
\pi^\bc \chi^{-1}(x_1 \connect \dots \connect x_n)
	= \pi^\bc \partial_\Omega(\Upsilon^{-1}(x_1)\udot\dots\udot\Upsilon^{-1}(x_k)).
\]
Let $y_i = \pi^\bc \Upsilon^{-1}(x_i)$.  Each $y_i$ has at least one
leg, since if
the $\partial_\Omega^{-1}$ of $\Upsilon^{-1}$ eats all the legs of
$\chi^{-1}x_i$, it also creates a vacuum diagram which is killed by
$\pi^\bc$.  Then
\[
\pi^\bc \partial_\Omega(y_1 \dots y_k) =
    \pi^\bc \langle \Omega_a, \Delta_{ab}(y_1 \dots y_k)\rangle_a.
\]
Let $\Delta_{ab} y_i = (y_i)_a + z_i$; diagrams in $z_1$ have at least
one $b$ leg.  We see that
\begin{align*}
\pi^\bc\langle\Omega_a, (y_1)_a \Delta_{ab}(y_2 \dots y_n)\rangle_a
 &= \pi^\bc\langle(\partial_{y_1}\Omega)_a, \Delta_{ab}(y_2\dots y_k)\rangle_a
	&&\text{by Lemma~\ref{lem:mult-diffop-dual}} \\
 &= 0.   &&\text{by Equation~\ref{eq:pseudo-linear-2}}
\end{align*}
Therefore
\begin{align*}
\pi^\bc \partial_\Omega(y_1 \dots y_k)
   &= \pi^\bc(\langle \Omega_a, (y_1)_a \Delta_{ab}(y_2 \dots y_k)\rangle_a
    + \langle \Omega_a, z_1 \Delta_{ab}(y_2 \dots y_k)\rangle_a) \\
   &= \pi^\bc\langle \Omega_a, z_1 \Delta_{ab}(y_2 \dots y_k)\rangle_a \\
   & = \cdots \\
   & = \pi^\bc\langle \Omega_a, z_1 \dots z_k \rangle_a.
\end{align*}
Each $z_i$ has at least one leg labelled $b$, so the product has at
least $k$ legs labelled $b$ which are the legs in the result.
\end{proof}

Consider expanding the exponential $\exp_{\connect}(\strutn{}{}/2n)$
in the LHS of Equation~\ref{eq:n00}.  In the term with $k$ struts,
there is a factor of $1/n^k$ from the factors in front of the struts.  On
the other hand, by Lemma~\ref{lem:wheels-lemma}, the 
product has at least $k$ legs, or $k+1$ if there is a non-trivial
contribution from $\nu$.
Since the overall power of $n$ is $n^{k - \text{\# legs}}$, when we
evaluate at $n=0$ the term $\nu$ does not contribute at all.
\[
\psi^{(n)}(\nu \connect \exp_{\connect}(\frac{1}{2n}\strutn{}{}))|_{n=0}
    = \psi^{(n)}(\exp_{\connect}(\frac{1}{2n}\strutn{}{}))|_{n=0}.
\]
Now we want to pick out the term from $(\strutn{}{})^{\connect k}$ with
exactly $k$ legs.  We can do this computation explicitly using the
wheeling map $\Upsilon$. Alternatively, the
result (which is $\nu$) must be a diagram of degree $k$ and with $k$
legs, hence of
mapping degree 0: $\nu = \nu_0$.  $\Omega \cap \strutv{}{}$ was shown in
Chapter~\ref{chap:wheeling} the be the part of $Z(\OpenHopf)$ of
mapping degree 0.  Dropping the central strand from $\OpenHopf$ leaves
an unknot, so $\Omega =\nu_0 = \nu$.
This completes the proof of Theorem~\ref{thm:wheels}.\qed

\begin{exercise} Do the computation suggested above.  Show that
\[
\exp_{\connect}(\frac{1}{2}\strutn{}{}) =
	\Omega \udot \exp_{\udot} (\frac{1}{2}\strutn{}{}).
\]
\end{exercise}
\begin{hint} Use Lemma~\ref{lem:inspired-guess}.
\end{hint}

\section {From the unknot to the Hopf link}\label{sec:Hopf}

We will now apply Theorems~\ref{thm:wheeling}
and~\ref{thm:wheels} to compute the invariant of the Hopf link.
A little attention is required in order to perform the operations in the
correct order.  We start by computing the Kontsevich integral of the
$+1$ framed unknot.
\begin{align*}
Z(\bigcircle^{+1}) &= \Omega \connect \exp_{\connect}(\frac 12\isolatedchord)\\
	&= \partial_\Omega\left(\partial_\Omega^{-1}(\Omega) \udot
		\exp_{\udot}(\partial_\Omega^{-1}(\strutn{}{}))\right)
		&& \text{by Theorem~\ref{thm:wheeling}} \\
	&= \pi^\bc\partial_\Omega\left(\Omega \udot \exp(\strutn{}{})\right).
		&& \text{by Equation~\ref{eq:pseudo-linear-2}}
\end{align*}

To pass to the Hopf link, we double $Z(\bigcircle^{+1})$.  The
following lemma tells us how $\partial_\Omega$ intereacts with
doubling.  We use $\hat D$ as an
alternate notation for $\partial_D$ so that we can use a
subscript to indicate which variable the differential operator acts on.
\begin{lemma}\label{lem:double-diffop} For $C, D \in \cB$ with $C$ strutless,
\[
\Delta_{xy} \hat C(D) = \hat C_x(\Delta_{xy} D) = \hat C_y(\Delta_{xy} D).
\]
\qed
\end{lemma}
If we want to apply $\partial_\Omega^{-1}$ to both components of the
Hopf link, we can compute $\partial_\Omega^{-2}(Z(\bigcircle^{+1}))$. We
make an inspired guess.

\begin{lemma}\label{lem:inspired-guess}
 $\pi^\bc\partial_\Omega(\exp\frac{1}{2}\strutn{}{}) =
	\Omega\udot\exp(\frac{1}{2}\strutn{}{})$.
\end{lemma}
\begin{proof}
\begin{align*}
\pi^\bc\partial_\Omega(\exp(\frac12\strutn{}{}))
   &= \pi^\bc \langle\Omega_y, \exp(\frac{1}{2} \strutn{x+y}{x+y})\rangle_y
\\
   &= \pi^\bc\langle\Omega_y,
	\exp(\frac12\strutn xx)\exp(\strutn yx)\exp(\frac12\strutn yy)\rangle_y \\
   &= \pi^\bc\langle \partial_{\exp(\frac12\strutn{}{})}(\Omega)_y,
	\exp(\strutn xy)\rangle_y\udot\exp(\frac12\strutn xx)
			&& \text{by Lemma~\ref{lem:mult-diffop-dual}} \\
   &= \pi^\bc\langle\Omega_y, \exp(\strutn yx)\rangle_y\udot\exp(\frac12\strutn xx)
			&& \text{by Equation~\ref{eq:pseudo-linear-2}} \\
   &= \Omega\udot\exp(\frac12\strutn{}{}).
\end{align*}
\end{proof}

As a corollary, we see that
\begin{equation} \label{eq:inspired-guess}
\pi^\bc\partial_\Omega^{-2}(Z(\bigcircle^{+1})) = \exp(\frac12\strutn{}{}).
\end{equation}
We now compute.
\begin{align*}
\pi^\bc \Delta_{xy}(\hat\Omega^{-2}Z(\bigcircle^{+2}))
   &= \pi^\bc \hat\Omega_x^{-1}\hat\Omega_y^{-1} Z(\,{}_x^{+1}\!\Hopf_y^{+1})
	&&\text{by Lemma~\ref{lem:double-diffop} and Proposition~\ref{prop:disconn-cable}} \\
   &= \pi^\bc\Delta_{xy}(\exp(\frac12\strutn{}{}))
	&&\text{by Equation~\ref{eq:inspired-guess}} \\
   &= \exp(\strutn{x}{y})\exp(\frac12\strutn xx)\exp(\frac12\strutn yy).
\end{align*}
Apply $\Upsilon_x \circ \Upsilon_y$ to both sides.  We see that
\begin{align*}
Z(\,{}_x^{+1}\!\Hopf_y^{+1}) &= \pi^\bc Z(\,{}_x^{+1}\!\Hopf_y^{+1}) \\
	&= \pi^\bc\Upsilon_x \circ \Upsilon_y(\exp(\strutn{x}{y})\udot\exp(\frac12\strutn xx)
			\udot\exp(\frac12\strutn yy)) \\
	&= \pi^\bc \Upsilon_x \circ \Upsilon_y(\exp(\strutn{x}{y}))
		\connect \exp_{\connect}(\frac12\isolatedchord_x)
		\connect \exp_{\connect}(\frac12\isolatedchord_y) \\
\intertext{so}
Z(\hopf{x}{y}) &= \pi^\bc \Upsilon_x \circ \Upsilon_y(\exp(\strutn{x}{y})).
\end{align*}
This is the first equality of Theorem~\ref{thm:hopf-link}.
For the second equality,
\[
\hat\Omega_y(\exp(\strutn xy)) = \Omega_x \udot \exp(\strutn xy).
\]
so
\[
Z(\hopf{x}{y}) = \pi^\bc \Upsilon_x(\exp(\strutn{y}{x})\Omega_x).
\]
For the last equality of the theorem, multipicativity of
$\Upsilon$ implies that
\begin{align*}
\pi^\bc\Upsilon_x(\exp(\strutn yx)\Omega_x) 
  &= \pi^\bc(\Upsilon_x(\exp(\strutn yx)) \connect \Upsilon_x(\Omega_x)) \\
  &= \pi^\bc(\Upsilon_x(\exp(\strutn yx))) \connect \chi(\Omega_x)
  &= \chi(\exp(\strutn yx) \udot \Omega_y) \connect \chi(\Omega_x).
\end{align*}
But by Equation~\ref{eq:11tangle},
\[
Z(\openhopf{x}{y}) = Z(\hopf{x}{y}) \connect \Omega_x^{-1}
	= \exp(\strutn yx) \udot \Omega_y.
\]
This completes the proof of Theorem~\ref{thm:hopf-link}.\qed
\appendix
\chapter{Cyclic invariance of vertices}\label{app:cyclic}

In this appendix, we tie up a loose end from
Section~\ref{sec:weight-systems}.  In that section, the diagrams we
found naturally from the finite-type theory were Jacobi diagrams with
an additional structure, a ``routing''.  Recall that a routing of a
Jacobi diagram is a choice of two edges incident to each internal
trivalent verex.

Here we will find some good routings (routings that
corresponds to $n$-circuses) for any boundary-connected Jacobi
diagram and show that all such routings are equal
modulo the Antisymmetry and Splitting relations.

First we will give a convenient criterion for a routing to be good in
the above sense.

\begin{lemma}\label{lem:pull-apart}
A routing of a Jacobi diagram $D$ is good if there exists an
ordering of the internal vertices of $D$ so that every vertex has at
least one neighbor along a distinguished edge which is either an
external vertex or comes earlier in the ordering.
\end{lemma}

\begin{proof}
We need to check that after resolving each vertex like
\[
\figcent{vertex}{0.5} \rightsquigarrow
	\figcent{vertex-expand}{0.5}.
\]
we are left with a collection of $n$ double lassos that are trivial in
$S^3$ once you forget the knot.  The internal vertices can be pulled
apart one by one in the specified order: at each stage, at least one
side of the encircled lasso ends in a loop with nothing inside.
\end{proof}

In fact, the condition in the lemma is an if and only if, but we do
not need that.

\begin{definition} With respect to an ordering of the internal
vertices of a Jacobi diagram, a vertex $v_1$ is {\em younger} than a
vertex $v_2$ if $v_1$ is an internal vertex or $v_1$ and $v_2$ are
both internal and $v_1$ is earlier than $v_2$ in the ordering.  An
ordering of the internal vertices is {\em good} if every vertex has a
younger neighbor.  A routing and an ordering of the internal vertices
are {\em compatible} if they satisfy the condition of
Lemma~\ref{lem:pull-apart}, i.e., if one of the two distinguished
neighbors of each vertex is younger.
\end{definition}

\begin{proposition}
Every boundary connected Jacobi diagram $D$ has at least one good routing.
\end{proposition}

\begin{proof} Because $D$ is boundary connected, there is a good ordering
on the internal vertices.  (Order the vertices from the
external vertices on the knot inward.)  Every good
ordering has a comptabile routing.
\end{proof}

\begin{lemma}\label{lem:good-order}
Any two good orderings are related by a series of transpositions of
vertices adjacent in the ordering.
\end{lemma}

\begin{proof} The minimal vertex in the first ordering must have an
external vertex as a neighbor.  Therefore, there is no obstruction to
moving this vertex to the first position in the second ordering by a
series of transpositions.  We can repeat this for each vertex in turn.
\end{proof}

\begin{lemma}\label{lem:good-route}
Any two two routings compatible with the same ordering are
equivalent modulo the vertex and antisymmetry relations.
\end{lemma}

\begin{proof}
We need to check that the two possibilities for the routing at each
internal vertex are equivalent.  This follows from two applications of
the vertex relation:
\[
\figcent{cyclic-trivalent-1}{0.5} =
\figcent{cyclic-trivalent-2}{0.5} -
\figcent{cyclic-trivalent-3}{0.5} =
\figcent{cyclic-trivalent-4}{0.5}
\]
where $i$ is younger than $j$ and the order of the two vertices in the
middle sum is chosen depending on which neighbor of $i$ is
younger, or is irrelevant if $i$ is an external vertex.
\end{proof}

\begin{proposition}
All routings of a Jacobi diagram $D$ compatible with any ordering are
equal modulo the antisymmetry and vertex relations.
\end{proposition}

\begin{proof} Pick a good ordering for each routing.  We can adjust
the routing while keeping the ordering fixed by
Lemma~\ref{lem:good-route}, so we just need to check that we can
change the ordering.  The two
orderings are related by a chain of good orderings related by adjacent
transpositions by Lemma~\ref{lem:good-order}.  At each transposition,
the two vertices involved must have younger neighbors which are not
each other (otherwise we could not exchange them); a routing in which
these two neighbors are on distinguished edges is compatible with
both orderings.
\end{proof}

\begin{remark} Goussarov~\cite{Goussarov:KnottedGraphs} and
Habiro~\cite{Habiro:Claspers} independently found a topological theory
of Jacobi diagrams in which the vertices are cyclically invariant from
the beginning.
\end{remark}

\chapter{Fixing the coefficients in $\Omega$}\label{app:coefficients}
To fix the coefficients in $\Omega$ (in terms of
Chapter~\ref{chap:wheeling}, to show that $\Omega' = \Omega$), we look
at
\[
\partial_{\Omega'}(\exp(\frac12 \strutn{}{}))
\]
as we did in Chapter~\ref{chap:wheels}.  In that chapter we projected
out the vacuum diagrams from the result, while here we will look at
the part that is only vacuum diagrams.  In particular, we have 
the following equality, which we will call ``Sawon's
identity~\cite{HitchinSawon:CurvatureHyperKahler}'':

\begin{equation}
\label{eq:sawon}
\langle \Omega', (\strutn{}{})^n \rangle = (\frac{1}{24}\ThetaGraph)^n.
\end{equation}

\begin{proof}
Proceed by induction on $n$.  The result is trivial for $n=0$.
\begin{align*}
\langle \Omega', (\strutn{}{})^n\rangle 
   &= \langle \Omega', \strutn{}{} \udot (\strutn{}{})^{n-1}\rangle \\
   &= \langle \partial_{\strutu{}{}} (\Omega'), (\strutn{}{})^{n-1}\rangle
	&&\text{by Lemma~\ref{lem:mult-diffop-dual}} \\
   &= \frac{1}{24}\ThetaGraph\,\langle\Omega', (\strutn{}{})^{n-1}\rangle
	&&\text{by Lemma~\ref{lem:pseudo-linear} and explicit computation} \\
   &= \left(\frac{1}{24}\ThetaGraph\right)^n
	&&\text{by induction}
\end{align*}
\end{proof}

Equation~\ref{eq:sawon} is already enough to fix the
coefficients in $\Omega'$. 

\begin{lemma}\label{lem:sl2rel}
In the Lie algebra $\sltwo$, we have the following relations:
\begin{equation*}
\bigcirc \equiv 3 \qquad
\IGraph \equiv \smoothing - \crossing
\end{equation*}
\end{lemma}
\begin{proof}
The first relation says that $\sltwo$ is 3-dimensional.  For the
second relation, note that both sides, considered as elements of
$\End(\sltwo \otimes \sltwo)$, are multiples of the projection onto
the antisymmetric part $\sltwo \wedge \sltwo$ (which is
3-dimensional).  A little computation fixes the constant.  (Note that
the constant depends on the metric.  Here we use $\langle x,y\rangle =
-\tr(xy)$, where the trace is taken in the adjoint representation.)
\end{proof}

Apply the $\sltwo$ relations above to both sides of Sawon's identity.
For the RHS, we find that that $\ThetaGraph \equiv 6$.  For
the LHS, we will use the following two lemmas.

\begin{lemma}\label{lem:sl2-wheel}
Modulo the $\sltwo$ relations, $\omega_{2n} \equiv 2(\strutn{}{})^n$.
\end{lemma}
\begin{proof}
Proceed by induction.  This is a straightforward computation for
$n=1$.  For $n>1$, compute as follows:
\begin{align*}
\omega_{2n} = \figcent{sl2-wheel-1}{0.12}
 	  = \figcent{sl2-wheel-2}{0.12} - \figcent{sl2-wheel-3}{0.12}
	  = \strutn{}{} \udot \figcent{sl2-wheel-4}{0.12}
	  = \strutn{}{} \udot\omega_{2n-2}.
\end{align*}
\end{proof}

\begin{lemma}\label{lem:sl2-gluing}
Modulo the $\sltwo$ relations, $\langle (\strutn{}{})^n, (\strutn{}{})^n \rangle = (2n+1)!$.
\end{lemma}
\begin{proof}
Proceed by induction.  The statement is trivial for $n=0$.  For $n>0$,
the two ends of the first strut on the left hand side can either
connect to the two ends of a single right hand strut or they can connect to
two different struts.  These happen in $2n$ and $2n\cdot(2n-2)$ ways,
respectively.  (Note that there are $2n\cdot(2n-1)$ ways in all of
gluing these two legs.)  We therefore have
\[
\figcent{sl2-glue-1}{0.15} = 2n\cdot\figcent{sl2-glue-2}{0.15}
	+ 2n \cdot (2n-2) \cdot \figcent{sl2-glue-3}{0.15}
\]
and
\begin{align*}
\bigl\langle(\Strutn)^n, (\Strutn)^n\bigr\rangle
  &= \bigl(2n\bigcircle + 2n\cdot(2n-2)\bigr)
	\bigl\langle(\Strutn^{n-1},(\Strutn)^{n-1}\bigr\rangle \\
  &\equiv 2n\cdot(2n+1)\bigl\langle(\Strutn^{n-1},(\Strutn)^{n-1}\bigr\rangle \\
  &\equiv (2n+1)! &&\text{by induction.}
\end{align*}
\end{proof}

\begin{proposition}\label{prop:coefficients}
$\sum a_n x^n = \frac{1}{2}\log \frac{\sinh (x/2)}{x/2}$.
\end{proposition}

\begin{proof}
By Lemma~\ref{lem:sl2-wheel}, we find
\[
\Omega' = \exp(\sum_n a_{2n} \omega_{2n})
  \equiv \exp(\sum_n 2a_{2n}(\Strutn)^n).
\]
Set $f(x) = \exp(2\sum a_{2n} x^n) = \sum f_n x^n$.  Then by
Lemma~\ref{lem:sl2-gluing},
\begin{align*}
\langle\Omega',(\Strutn)^n\rangle &\equiv
	 \langle f(\Strutn), (\Strutn)^n\rangle
	= \langle f_n(\Strutn)^n, (\Strutn)^n\rangle
	\equiv f_n (2n+1)! \\
  &= \left(\frac{1}{24}\ThetaGraph\right)^n \equiv \frac{1}{4^n}.
\end{align*}
so
\begin{align*}
f_n &= \frac{1}{4^n (2n+1)!} \\
f(x) &= \frac{\sinh(\sqrt{x}/2)}{\sqrt{x}/2} \\
\exp\left(2\sum_n a_n x^n\right) &= \frac{\sinh (x/2)}{x/2} \\
\sum_n a_n x^n &= \frac{1}{2}\log\frac{\sinh(x/2)}{x/2}.
\end{align*}
\end{proof}
\bibliographystyle{hamsplain}
\bibliography{knots}

\providecommand{\bysame}{\leavevmode\hbox to3em{\hrulefill}\thinspace}
\begin{thebibliography}{10}

\bibitem{AlekseevMeinrenken:NonCommutativeWeil}
Anton Alekseev and Eckhard Meinrenken, \emph{The non-commutative {W}eil
  algebra}, Invent. Math. \textbf{139} (2000), no.~1, 135--172,
  \mbox{arXiv:math.DG/9903052}.

\bibitem{ADS:KontsevichQuantization}
Martin Andler, Alexander Dvorsky, and Siddhartha Sahi, \emph{{K}ontsevich
  quantization and invariant distributions on {L}ie groups},
  \mbox{arXiv:math.QA/9910104}.

\bibitem{BarNatan:OnVassiliev}
Dror Bar-Natan, \emph{On the {V}assiliev knot invariants}, Topology \textbf{34}
  (1995), 423--472.

\bibitem{BGRT:WheelsWheeling}
Dror Bar-Natan, Stavros Garoufalidis, Lev Rozansky, and Dylan~P. Thurston,
  \emph{Wheels, wheeling, and the {K}ontsevich integral of the unknot},
  \mbox{arXiv:q-alg/9703025}.

\bibitem{Bar-NatanLawrence:RationalSurgery}
Dror Bar-Natan and Ruth Lawrence, \emph{A rational surgery formula for the
  {LMO} invariant}, in preparation.

\bibitem{Chen:Iterated}
Kuo-tsai Chen, \emph{Iterated integrals of differential forms and loop space
  homology}, Ann. of Math. (2) \textbf{97} (1973), 217--246.

\bibitem{Deligne:letter}
Pierre Deligne, letter to D. Bar-Natan, January 1996.

\bibitem{Duflo:Operateurs}
Michel Duflo, \emph{Opérateurs différentiels bi-invariants sur un groupe de
  {L}ie}, Ann. scient. École Norm. Sup. \textbf{10} (1977), 265--288.

\bibitem{Goussarov:KnottedGraphs}
Mikhail Goussarov, \emph{Knotted graphs and a geometrical technique of
  $n$-equivalence}, Tech. report, POMI Sankt Petersburg, 1995.

\bibitem{Goussarov:3Manifolds}
Mikhail Goussarov, \emph{Finite type invariants and $n$-equivalence of
  $3$-manifolds}, C. R. Acad. Sci. Paris S\'er. I Math. \textbf{329} (1999),
  no.~6, 517--522.

\bibitem{Habiro:Claspers}
Kazuo Habiro, \emph{Claspers and finite type invariants of links}, Geom. Topol.
  \textbf{4} (2000), 1--83.

\bibitem{HitchinSawon:CurvatureHyperKahler}
Nigel Hitchin and Justin Sawon, \emph{Curvature and characteristic numbers of
  hyperkähler manifolds}, August 1999, \mbox{arXiv:math.DG/9908114}.

\bibitem{Kontsevich:DeformationQuantization}
Maxim Kontsevich, \emph{{Deformation quantization of Poisson manifolds, I}},
  \mbox{arXiv:q-alg/9709040}.

\bibitem{Kuperberg:Invertibility}
Greg Kuperberg, \emph{{Detecting knot invertibility}}, J. Knot Theory
  Ramifications \textbf{5} (1996), 173--181, \mbox{arXiv:q-alg/9712048}.

\bibitem{LeMurakami:Universal}
Thang T.~Q. Le and Jun Murakami, \emph{The universal {V}assiliev-{K}ontsevich
  invariant for framed oriented links}, preprint, Max-Planck-Institut Bonn,
  December 1993.

\bibitem{LeMurakami:Parallel}
\bysame, \emph{Parallel version of the universal {K}ontsevich-{V}assiliev
  invariant}, J. Pure Alg. Appli. \textbf{212} (1997), 271--291.

\bibitem{Lescop:Kontsevich}
Christine Lescop, \emph{Kontsevich integral}, lecture notes, {Grenoble Summer
  School}, June 1999.

\bibitem{MO:TQFTUniversal}
Jun Murakami and Tomatada Ohtsuki, \emph{Topological quantum field theory for
  the universal quantum invariant}, Commun. Math. Phs. \textbf{188} (1997),
  501--520.

\bibitem{Thurston:Torus}
Dylan Thurston, \emph{Torus actions for the {LMO} invariant}, in preparation.

\bibitem{Vogel:Structures}
Pierre Vogel, \emph{Algebraic structures on modules of diagrams}, Tech. report,
  Université Paris VII, July 1995.

\end{thebibliography}
\end{document}